\documentclass[11pt]{article}
\usepackage[margin=2.3cm]{geometry} 

\usepackage[utf8]{inputenc}
\usepackage{amsfonts}
\usepackage{amssymb}
\usepackage{tikz-cd} 
\usepackage{graphicx}
\usepackage{mathtools} 
\usepackage{setspace} 
\usepackage{colonequals} 
\usepackage{anyfontsize} 
\usepackage{quiver}
\usepackage{physics}
\usepackage{tikz}
\usepackage{mathdots}
\usepackage{yhmath}
\usepackage{cancel}
\usepackage{array}
\usepackage{multirow}
\usepackage{amssymb}
\usepackage{tabularx}
\usepackage{extarrows}
\usepackage{booktabs}
\usetikzlibrary{fadings}
\usetikzlibrary{patterns}
\usetikzlibrary{shadows.blur}
\usetikzlibrary{shapes}
\usepackage{color}

\usepackage{microtype} 
\usepackage[nottoc]{tocbibind} 
\usepackage[none]{hyphenat} 
\usepackage{subcaption} 

\usepackage{titlesec}
\titleformat*{\section}{\Large\bfseries\sffamily} 
\titleformat*{\subsection}{\bfseries\sffamily}
\titleformat*{\subsubsection}{\bfseries\sffamily}

\usepackage{hyperref}
\usepackage{setspace}

\usepackage{hyperref}
\hypersetup{
    colorlinks=true,
    allcolors=blue
}

\usepackage{lmodern}

\usepackage{lscape}
\usepackage{graphicx}
\usepackage{amssymb}
\usepackage{mathtools}
\usepackage{amsmath}
\usepackage{multirow}
\usepackage{amsthm}
\usepackage{epstopdf}
\usepackage{youngtab}
\usepackage{bm}
\usepackage{tikz}
\usetikzlibrary{cd}
\usetikzlibrary{calc}
\usetikzlibrary{arrows}
\usetikzlibrary{arrows.meta}
\usetikzlibrary{shapes.arrows}
\usetikzlibrary{decorations.pathmorphing}
\usepackage{url}
\usepackage{array}
\begingroup
    \makeatletter
    \@for\theoremstyle:=definition,remark,plain\do{%
        \expandafter\g@addto@macro\csname th@\theoremstyle\endcsname{%
            \addtolength\thm@preskip\parskip
            }%
        }
\endgroup

\allowdisplaybreaks

\def\dashmapsto{\mathrel{\mapstochar\dashrightarrow}} 

\newcommand{\dash}{\dashrightarrow}
\newcommand{\nin}{\noindent}

\newtheorem{thm}{Theorem}[section]
\newtheorem{prop}[thm]{Proposition}
\newtheorem{propdefn}[thm]{Proposition--Definition}
\newtheorem{lem}[thm]{Lemma}
\newtheorem{cor}[thm]{Corollary}

\theoremstyle{definition}
\newtheorem{defn}[thm]{Definition}
\newtheorem{ex}[thm]{Example}
\newtheorem{rmk}[thm]{Remark}

\newtheorem*{nota}{Notation}

\numberwithin{equation}{section}

\newcommand{\p}{\mathbb{P}}
\newcommand{\q}{\mathbb{Q}}
\newcommand{\C}{\mathbb{C}}
\newcommand{\F}{\mathbb{F}}

\newcommand{\ZZ}{\mathbb{Z}}

\newcommand{\bF}{\mathbf{F}}

\newcommand{\oo}{\mathcal{O}}

\newcommand{\h}{\mathcal{H}}
\newcommand{\n}{\mathcal{N}}

\DeclareMathOperator{\Bl}{Bl}

\DeclareMathOperator{\Spec}{Spec}

\DeclareMathOperator{\Id}{Id}
\DeclareMathOperator{\Exc}{Exc}

\DeclareMathOperator{\PGL}{PGL}

\DeclareMathOperator{\Bir}{Bir}

\DeclareMathOperator{\vp}{vp}
\DeclareMathOperator{\Bs}{Bs}
\DeclareMathOperator{\uBs}{\underline{Bs}}
\DeclareMathOperator{\Aut}{Aut}

\DeclareMathOperator{\divv}{div}

\DeclareMathOperator{\Pic}{Pic}

\DeclareMathOperator{\pr}{pr}
\DeclareMathOperator{\mult}{mult}
\DeclareMathOperator{\Dec}{Dec}
\DeclareMathOperator{\Ine}{Ine}

\DeclareMathOperator{\Sing}{Sing}

\DeclareMathOperator{\Dom}{Dom}

\DeclareMathOperator{\rk}{rank}

\title{\textbf{\textsf{On the decomposition group of a nonsingular plane cubic by a log Calabi-Yau geometrical perspective}}}
\author{Eduardo Alves da Silva
\thanks{Universität Basel, Basel - Switzerland. \\
\hspace{1.0cm} \textit{Email}: \texttt{eduardo.alvesdasilva.math@gmail.com}}}
\date{} 

\begin{document}

\maketitle

\renewcommand*\abstractname{\textsf{Abstract}}

\begin{abstract}
This paper aims to study the decomposition group of a nonsingular plane cubic under the light of the log Calabi-Yau geometry.  Using this approach we prove that an appropriate algorithm of the Sarkisov program in dimension 2 applied to an element of this group is automatically volume preserving. From this, we deduce some properties of the (volume preserving) Sarkisov factorization of its elements. We also confirm an expectation by Blanc, Pan and Vust on the splitting property of the canonical complex of a nonsingular plane cubic. Within a similar context in dimension 3, we exhibit in detail an interesting counterexample for a possible generalization of a theorem by Pan in which there exists a Sarkisov factorization obtained algorithmically that is not volume preserving. 

\end{abstract}




{
  \hypersetup{linkcolor=black}

\renewcommand*\contentsname{\textsf{Contents}}

\setcounter{tocdepth}{1}

\renewcommand{\baselinestretch}{0}

\textsf{\tableofcontents}

\renewcommand{\baselinestretch}{1.0}\normalsize
  
}

\section{Introduction}

In the context of algebraic geometry, decomposition and inertia groups are special subgroups of the Cremona group that preserve a certain subvariety of $\mathbb{P}^n$ as a set and pointwise, respectively. In \cite{cas,mm} and more recently in \cite{giz,pan,bpv1,bpv2,bl1,bl2,hz,dhz,piñ}, there exist numerous interesting results and descriptions of these groups in several cases. 

In the particular case where this fixed subvariety is a hypersurface $D_{n+1} \subset \p^n$ of degree $n+1$, we have that $(\p^n,D_{n+1})$ is an example of a \textit{Calabi-Yau pair}, that is, a pair $(X, D)$ with mild singularities consisting of a normal projective variety $X$ and a reduced Weil divisor on $X$ such that $K_X + D \sim 0$. In other words, regarding $n= \dim(X)$, there exists a meromorphic volume form $\omega = \omega_{X,D} \in \Omega^n_{X}$ up to nonzero scaling, such that $D+ \divv(\omega)=0$. 

Under restrictions on the singularities of $(\p^n,D)$, the decomposition group of the hypersurface $D$, denoted by $\Bir(\p^n,D)$ or simply $\Dec(D)$, coincides with the group of birational self-maps of $\p^n$ that preserve the volume form $\omega$, up to nonzero scaling. See Proposition \ref{vp can pair}. Such maps are naturally called \textit{volume preserving}.

This notion of Calabi-Yau pair allows us to use new tools to deal with the study of these groups and (re)interpret some results as statements about the birational geometry of the pair. One of these tools is the so-called volume preserving Sarkisov program, a result by Corti \& Kaloghiros \cite[Theorem~1.1]{ck} valid in all dimensions. From the perspective of log Calabi-Yau geometry, this result is a generalization of the standard Sarkisov program \cite{cor1,hm1} with additional structures.

The Sarkisov program asserts that we can decompose any birational map $\phi \in \Bir(\p^n)$ into a composition of finitely many elementary links between Mori fibered spaces, the so-called \textit{Sarkisov links}:

\begin{center}
\begin{tikzcd}
	{\mathbb{P}^n=X_0} & {X_1} & \cdots & {X_{m-1}} & {X_m=\mathbb{P}^n} \\
	{\text{Spec}(\mathbb{C})=Y_0} & {Y_1} && {Y_{m-1}} & {Y_m=\text{Spec}(\mathbb{C})}
	\arrow["{\phi_1}", dashed, from=1-1, to=1-2]
	\arrow["{\phi_2}", dashed, from=1-2, to=1-3]
	\arrow["{\phi_{m-1}}", dashed, from=1-3, to=1-4]
	\arrow["{\phi_m}", dashed, from=1-4, to=1-5]
	\arrow[from=1-1, to=2-1]
	\arrow[from=1-2, to=2-2]
	\arrow[from=1-4, to=2-4]
	\arrow[from=1-5, to=2-5]
	\arrow["\phi", curve={height=-24pt}, dashed, from=1-1, to=1-5]
\end{tikzcd}.
\end{center}

The volume preserving version of this result established in \cite{ck} states that, if the birational map $\phi$ is volume preserving for certain log canonical Calabi-Yau pairs $(\p^n,D)$ and $(\p^n,D')$, then there exists a factorization as above, and divisors $D_i \subset X_i$ making each $(X_i,D_i)$ a mildly singular Calabi-Yau pair, and each $\phi_i$ a volume preserving map between Calabi-Yau pairs. See Definition \ref{defn crep bir map}.

Our main result is the following:

\begin{thm}[See Theorem \ref{thm sark=vp cubic}]
Let $C \subset \p^2$ be a nonsingular cubic. The standard Sarkisov program applied to an element of $\Dec(C)$ is automatically volume preserving.
\end{thm}

The main fact used in the proof of this result is Theorem \ref{pan thm 1.3} due to Pan \cite{pan}, which asserts that the base locus of an element $\phi \in \Dec(C) \setminus \PGL(3,\C)$ is contained in the curve $C$. Furthermore, by employing the volume preserving variant of the Sarkisov program, it becomes feasible to demonstrate a broader fact: all the infinitely near base points of an element of $\Dec(C)$ belong to the strict transforms of $C$. See Lemma \ref{lem inf base locus cont in C}. This will guarantee that when running the (volume preserving) Sarkisov program, all the surfaces involved, together with the corresponding strict transforms of $C$, are always Calabi-Yau pairs. 

By \cite[Theorem 6]{giz} and \cite[Theorem 2.2]{og}, we have a natural exact sequence induced by the natural action $\rho$ of $\Dec(C)$ on $C$ 

\begin{equation}\label{cc int}
    1 \longrightarrow \Ine(C) \longrightarrow \Dec(C) \overset{\rho} \longrightarrow \Bir(C) \longrightarrow 1 ,
\end{equation}

\noindent where the inertia group $\Ine(C)$ is defined as $\ker(\rho)$. Blanc, Pan \& Vust \cite{bpv1} remarked about the splitting property of this sequence. Notice that in particular, $C$ is an elliptic curve and therefore also an algebraic group. One has $\Bir(C)=\Aut(C)=C \rtimes \mathbb{Z}_d$, where $C$ is the group of translations and $d \in \{2,4,6\}$, depending on the $j$ invariant of $C$. In the proof of \cite[Theorem 2.2]{og}, the surjectivity of $\rho$ is proved by exhibiting a set-theoretical section $C \hookrightarrow  \Dec(C)$. We show that this set-theoretical section, however, is not a group homomorphism and hence it is not a partial splitting of \ref{cc int}. From this, we are able to produce a lot of elements in $\Ine(C)$.

More generally, we show the following which confirms the expectation in \cite{bpv1}: 

\begin{thm}[See Theorem \ref{non split can comp}]
    The canonical complex \ref{cc int} of the pair $(\p^2,C)$ does not admit any splitting as groups at $C$ when we write $\Aut(C) = C \rtimes \mathbb{Z}_d$.
    
\end{thm}

Within a similar context in higher dimension, it is natural to ask ourselves about a generalization of the Theorem \ref{pan thm 1.3} and of the Theorem \ref{thm sark=vp cubic}. In dimension $3$ we exhibit in detail an interesting counterexample for these both questions arising from the decomposition group of a general quartic surface with a single canonical singularity of type $A_1$. See Section \ref{vp sark x sark fac}.

\subsection{Structure of the paper}

Throughout this paper, our ground field will be $\C$, or more generally, any algebraically closed field of characteristic zero. Concerning general aspects of birational geometry, we refer the reader to \cite{km}. 

In Section \ref{log CY geometry and vp Sark} we will give an overview of Calabi-Yau pairs and the Sarkisov program and we will introduce some natural classes of singularities of pairs. In Section \ref{dec group cubic}, we will approach the decomposition group of a nonsingular plane cubic via log Calabi-Yau geometry aiming the proof the Theorem \ref{thm sark=vp cubic} and its corollaries. We will also deal with the open question left in \cite{bpv1}.
In Section \ref{vp sark x sark fac}, we will explain in detail the counterexample for a generalization of Theorem \ref{pan thm 1.3} comparing its possible Sarkisov factorizations, volume preserving or not.

\subsection{Acknowledgements}

I would like to thank my PhD advisor Carolina Araujo for guiding me in the process of construction of this paper and for her tremendous patience in our several discussions. I also thank CNPq (Conselho Nacional de Desenvolvimento Científico e Tecnológico) for the financial support during my studies at IMPA, and FAPERJ (Fundação de Amparo à Pesquisa do Estado do Rio de Janeiro) for supporting my attendance in events that were very important to finding people to discuss this work.

\section{Log Calabi-Yau geometry and the volume preserving Sarkisov program}\label{log CY geometry and vp Sark}

In this section, we give an overview of the concepts and ideas involving the geometry of log Calabi-Yau pairs and the volume preserving Sarkisov program.

\begin{defn}
A \textit{log Calabi-Yau pair} is a log canonical pair $(X,D)$ consisting of a normal projective variety $X$ and a reduced Weil divisor $D$ on $X$ such that $K_X + D \sim 0$. This condition implies the existence of a top degree rational differential form $\omega = \omega_{X,D} \in \Omega_X^n$, unique up to nonzero scaling, such that $D + \divv(\omega)=0$. By abuse of language, we call this differential the \textit{volume form}.
\end{defn}

From now on, we will call a log Calabi-Yau pair simply a Calabi-Yau pair. Sometimes in a more general context, it is admitted that $D$ is a $\q$-divisor and that $K_X + D \sim_{\q} 0$, but we will not need this generality here. Since $K_X + D \sim 0$, we have that $K_X + D$ is readily Cartier,
and hence all the discrepancies with respect to the pair $(X,D)$ are integer numbers.

We now introduce some classes of pairs taking into account the singularities of the ambient varieties and the divisors.

\begin{defn}\label{defn pairs}
    We say that a pair $(X,D)$ is $(t,c)$, respectively, $(t,lc)$, if $X$ has terminal singularities and the pair $(X,D)$ has canonical, respectively log canonical singularities. We say that a pair $(X,D)$ is $\q$-factorial if $X$ is $\q$-factorial.  
\end{defn}

If $(X,D)$ is (t,lc), then $a(E,X,D_X) \leq 0$ implies $a(E,X,D_X) = -1$ or $0$.

\begin{propdefn}[cf. \cite{km} Lemma 2.30]\label{prop crep bir morp}
    A proper birational morphism $f \colon (Z,D_Z) \rightarrow (X,D_X)$ is called crepant if $f_*D_Z=D_X$ and $f^*(K_X + D_X) \sim K_Z+ D_Z$. The term ``crepant'' (coined by Reid) refers to the fact that for every divisor $E$ over $X$, one has $a(E,Z,D_Z)=a(E,X,D_X)$. 
\end{propdefn}

\begin{defn}\label{defn crep bir map} Let $(X,D_X) $ and $(Y,D_Y)$ be log canonical pairs. We say that they are \textit{crepant birational equivalent} if there exists a birational map of pairs $\phi \colon (X,D_X) \dash (Y,D_Y)$ that admits a resolution
    \begin{center}
\begin{tikzcd}
	& {(Z,D_Z)} \\
	{(X,D_X)} && {(Y,D_Y)}
	\arrow["\phi", dashed, from=2-1, to=2-3]
	\arrow["p"', from=1-2, to=2-1]
	\arrow["q", from=1-2, to=2-3]
\end{tikzcd}
    \end{center}

\noindent in such a way that $p$ and $q$ are crepant birational morphisms.

\end{defn}

This definition is equivalent to asking that $a(E,X,D_X)=a(E,Y,D_Y)$ for every valuation $E$ of $K(X) \simeq K(Y)$ as in item \ref{2 prop vp map} in Proposition \ref{prop vp map}.

For Calabi-Yau pairs, the notion of \textit{crepant birational equivalence} becomes \textit{volume preserving equivalence}, since $\phi^*\omega_{Y,D_Y} = \lambda \omega_{X,D_X}$, for some $\lambda \in \C^*$. In this case, we call such $\phi$ a \textit{volume preserving} map.

As a consequence of these equivalences, we have the following:

\begin{prop}[cf. \cite{ck} Remark 1.7]\label{prop vp map}
Let $(X,D_X)$ and $(Y,D_Y)$ be Calabi-Yau pairs and $\phi \colon X \dash Y$ an arbitrary birational map. The following conditions are equivalent: 

\begin{enumerate}
    \item The map $\phi \colon (X,D_X) \dash (Y,D_Y)$ is volume preserving.

    \item\label{2 prop vp map} For all geometric valuations $E$ with center on both $X$ and $Y$, the discrepancies of $E$ with respect to the pairs $(X,D_X)$ and $(Y,D_Y)$ are equal: $a(E,X,D_X) = a(E,Y,D_Y)$.

    \item Let 
\[\begin{tikzcd}
	& {(Z,D_Z)} \\
	{(X,D_X)} && {(Y,D_Y)}
	\arrow["\phi", dashed, from=2-1, to=2-3]
	\arrow["p"', from=1-2, to=2-1]
	\arrow["q", from=1-2, to=2-3]
\end{tikzcd}\]
    be a common log resolution of the pairs $(X,D_X)$ and $(Y,D_Y)$. The birational map $\phi$ induces an identification $\phi_* \colon \Omega^n_X \xrightarrow{\sim} \Omega^n_Y$, where $n = \dim(X)=\dim(Y)$. 
    By abuse of notation, we write \begin{center}
$p^*(K_X+D_X)=q^*(K_Y+D_Y)$
\end{center}

to mean that for all $\omega \in \Omega^n_X$, we have \begin{center}
$p^*(D_X+\divv(\omega))=q^*(D_Y+\divv(\phi_*(\omega)))$.
\end{center}

The condition is: for some (or equivalently for any) common log resolution as above, we have \begin{center}
$p^*(K_X+D_X)=q^*(K_Y+D_Y)$.
\end{center}
    
\end{enumerate}
    
\end{prop}

\begin{defn}\label{Bir^vp}  
Let $(X,D)$ be a Calabi-Yau pair. The \textit{group of volume preserving birational self-maps} of $(X,D)$ is denoted by  $\Bir^{\vp}(X,D) \subset \Bir(X)$, and consists of all birational self-maps of $(X,D)$ that preserve the volume form.  
\end{defn}

\begin{defn}
    A \textit{Mori fibered Calabi-Yau pair} is a $\q$-factorial (t,lc) Calabi-Yau pair $(X,D)$ with a Mori fibered space structure on $X$, that is, a morphism $f \colon X \rightarrow S$, to a lower dimensional normal variety $S$, such that $f_*\oo_X=\oo_S$ (i.e. a contraction), $-K_X$ is $f$-ample, and $\rho(X/S)=\rho(X) - \rho(S)=1$. 
\end{defn}

\paragraph{Sarkisov program.} Since we will deal most of the time with $2$-dimensional pairs, we will recall the Sarkisov program in this case. We point out that each of the following elementary links has an analog in higher dimension. See \cite{cor1,hm1}. It is well known that the minimal models of rational surfaces are $\p^2$ and the geometrically ruled surfaces $\F_n$ with $n \neq 1$, also known as Hirzebruch surfaces. See \cite[Theorem V.10]{bea}. Such surfaces $\F_n$ are defined to be the $\p^1$-bundle over $\p^1$ given as the projective bundle of the rank two vector bundle $\oo_{\p^1} \oplus \oo_{\p^1}(n)$. In particular, one can show that $\F_0 \simeq \p^1 \times \p^1$ and $\F_1$ is isomorphic to $\p^2$ blown up at a point.

Seen as Mori fibered spaces, all these surfaces carry a structure morphism as follows: $\p^2 \rightarrow \Spec(\C)$, the projections on the two factors $p_1 \colon \p^1 \times \p^1 \rightarrow \p^1$ and $p_2 \colon \p^1 \times \p^1 \rightarrow \p^1$, $\F_n \rightarrow \p^1$.

\paragraph{Geometry of the Hirzebruch surfaces $\F_n$.} The Picard group of the surface $\F_n$ is isomorphic to $\ZZ \cdot [F] \oplus \ZZ \cdot [E]$, where $F$ is a fiber of the structure morphism and
\begin{center}
    $E = \begin{cases} \hfil
\text{the negative section}, & \text{if $n \geq 1$}\\
\text{any section with self-intersection 0} , & \text{if $n=0$}
\end{cases}$.
\end{center}

To simplify the notation, sometimes we will identify divisors with their corresponding classes in the Picard group. By abuse of notation, we will denote $ \ZZ \cdot F \oplus \ZZ \cdot E$ as $\langle F,E \rangle$. Furthermore, for $n=0$ and also abusing notation, we will call any section with 0 self-intersection a ``negative section''.\\

In the following description, we dispose the varieties of same Picard rank at the same height. Every birational self-map of the projective plane or minimal Hirzebruch surfaces is a composition of the following elementary maps:

\begin{enumerate}
    \item A \textit{Sarkisov link of type I} is a commutative diagram:
    \begin{center}
\begin{tikzcd}
	& {\mathbb{F}_1} \\
	{\mathbb{P}^2} & {\mathbb{P}^1} \\
	{\text{Spec}(\mathbb{C})}
	\arrow[from=2-1, to=3-1]
	\arrow["\pi"', from=1-2, to=2-1]
	\arrow[from=1-2, to=2-2]
	\arrow[from=2-2, to=3-1]
\end{tikzcd}
    \end{center}

    \noindent where $\pi^{-1} \colon \p^2 \dash \F_1 $ is a point blowup. 

    \item A \textit{Sarkisov link of type II} is a commutative diagram: 
    \begin{center}
\begin{tikzcd}
	{\mathbb{F}_n} & {\mathbb{F}_{n\pm 1}} \\
	{\mathbb{P}^1} & {\mathbb{P}^1}
	\arrow[from=1-1, to=2-1]
	\arrow["{\alpha_P}", dashed, from=1-1, to=1-2]
	\arrow[from=1-2, to=2-2]
	\arrow[Rightarrow, no head, from=2-1, to=2-2]
\end{tikzcd}
    \end{center}

    \noindent This is an elementary transformation $\alpha_P \colon \F_n \dash \F_{n\pm 1}$, by which we mean the blowup of a point $P \in \F_n$, followed by the contraction of the strict transform of the fiber through $P$.

\begin{figure}[h]
\centering
\begin{subfigure}{0.48\textwidth}
\resizebox{1.0\textwidth}{!}{

\tikzset{every picture/.style={line width=0.75pt}} 

\begin{tikzpicture}[x=0.75pt,y=0.75pt,yscale=-1,xscale=1]

\draw   (55.5,256) -- (216.3,256) -- (216.3,405) -- (55.5,405) -- cycle ;
\draw [color={rgb, 255:red, 208; green, 2; blue, 27 }  ,draw opacity=1 ][fill={rgb, 255:red, 208; green, 2; blue, 27 }  ,fill opacity=1 ][line width=1.5]    (54.45,330.8) -- (217.35,330.2) ;
\draw [color={rgb, 255:red, 126; green, 211; blue, 33 }  ,draw opacity=1 ][line width=1.5]    (155.4,254.4) -- (155.4,404.9) ;
\draw  [fill={rgb, 255:red, 144; green, 19; blue, 254 }  ,fill opacity=1 ] (152,294.6) .. controls (152,292.83) and (153.43,291.4) .. (155.2,291.4) .. controls (156.97,291.4) and (158.4,292.83) .. (158.4,294.6) .. controls (158.4,296.37) and (156.97,297.8) .. (155.2,297.8) .. controls (153.43,297.8) and (152,296.37) .. (152,294.6) -- cycle ;
\draw   (403.5,256) -- (564.3,256) -- (564.3,405) -- (403.5,405) -- cycle ;
\draw [color={rgb, 255:red, 208; green, 2; blue, 27 }  ,draw opacity=1 ][fill={rgb, 255:red, 208; green, 2; blue, 27 }  ,fill opacity=1 ][line width=1.5]    (402.45,330.8) -- (565.35,330.2) ;
\draw [color={rgb, 255:red, 144; green, 19; blue, 254 }  ,draw opacity=1 ][fill={rgb, 255:red, 144; green, 19; blue, 254 }  ,fill opacity=1 ][line width=1.5]    (503.4,254.4) -- (503.4,404.9) ;
\draw  [fill={rgb, 255:red, 126; green, 211; blue, 33 }  ,fill opacity=1 ] (500.2,329.65) .. controls (500.2,327.88) and (501.63,326.45) .. (503.4,326.45) .. controls (505.17,326.45) and (506.6,327.88) .. (506.6,329.65) .. controls (506.6,331.42) and (505.17,332.85) .. (503.4,332.85) .. controls (501.63,332.85) and (500.2,331.42) .. (500.2,329.65) -- cycle ;
\draw   (223,31) -- (383.8,31) -- (383.8,180) -- (223,180) -- cycle ;
\draw [color={rgb, 255:red, 208; green, 2; blue, 27 }  ,draw opacity=1 ][fill={rgb, 255:red, 208; green, 2; blue, 27 }  ,fill opacity=1 ][line width=1.5]    (221.95,105.8) -- (384.85,105.2) ;
\draw [color={rgb, 255:red, 126; green, 211; blue, 33 }  ,draw opacity=1 ][line width=1.5]    (298.7,65.3) .. controls (328.7,110.8) and (338.7,128.3) .. (308.2,163.3) ;
\draw [color={rgb, 255:red, 144; green, 19; blue, 254 }  ,draw opacity=1 ][line width=1.5]    (314.2,38.3) .. controls (333.7,71.8) and (329.2,83.8) .. (300.2,97.8) ;
\draw    (204.2,141.5) -- (140.05,221.84) ;
\draw [shift={(138.8,223.4)}, rotate = 308.61] [color={rgb, 255:red, 0; green, 0; blue, 0 }  ][line width=0.75]    (10.93,-3.29) .. controls (6.95,-1.4) and (3.31,-0.3) .. (0,0) .. controls (3.31,0.3) and (6.95,1.4) .. (10.93,3.29)   ;
\draw    (423,136) -- (481.48,202.89) ;
\draw [shift={(482.8,204.4)}, rotate = 228.84] [color={rgb, 255:red, 0; green, 0; blue, 0 }  ][line width=0.75]    (10.93,-3.29) .. controls (6.95,-1.4) and (3.31,-0.3) .. (0,0) .. controls (3.31,0.3) and (6.95,1.4) .. (10.93,3.29)   ;
\draw  [dash pattern={on 4.5pt off 4.5pt}]  (258.8,336) -- (378.8,336) ;
\draw [shift={(380.8,336)}, rotate = 180] [color={rgb, 255:red, 0; green, 0; blue, 0 }  ][line width=0.75]    (10.93,-3.29) .. controls (6.95,-1.4) and (3.31,-0.3) .. (0,0) .. controls (3.31,0.3) and (6.95,1.4) .. (10.93,3.29)   ;

\draw (570.5,324.3) node [anchor=north west][inner sep=0.75pt]  [color={rgb, 255:red, 0; green, 0; blue, 0 }  ,opacity=1 ]  {$E'$};
\draw (407,335.3) node [anchor=north west][inner sep=0.75pt]  [color={rgb, 255:red, 0; green, 0; blue, 0 }  ,opacity=1 ]  {$-( n-1)$};
\draw (486.5,258.3) node [anchor=north west][inner sep=0.75pt]  [color={rgb, 255:red, 0; green, 0; blue, 0 }  ,opacity=1 ]  {$0$};
\draw (496.5,237.3) node [anchor=north west][inner sep=0.75pt]    {$F'$};
\draw (483.5,306.8) node [anchor=north west][inner sep=0.75pt]    {$Q$};
\draw (387,238.3) node [anchor=north west][inner sep=0.75pt]    {$\mathbb{F}_{n-1}$};
\draw (222,324.3) node [anchor=north west][inner sep=0.75pt]  [color={rgb, 255:red, 0; green, 0; blue, 0 }  ,opacity=1 ]  {$E$};
\draw (59,331.8) node [anchor=north west][inner sep=0.75pt]  [color={rgb, 255:red, 0; green, 0; blue, 0 }  ,opacity=1 ]  {$-n$};
\draw (139.5,258.3) node [anchor=north west][inner sep=0.75pt]  [color={rgb, 255:red, 0; green, 0; blue, 0 }  ,opacity=1 ]  {$0$};
\draw (148,237.3) node [anchor=north west][inner sep=0.75pt]    {$F$};
\draw (135,284.8) node [anchor=north west][inner sep=0.75pt]    {$P$};
\draw (41,235.8) node [anchor=north west][inner sep=0.75pt]    {$\mathbb{F}_{n}$};
\draw (388,90.3) node [anchor=north west][inner sep=0.75pt]  [color={rgb, 255:red, 0; green, 0; blue, 0 }  ,opacity=1 ]  {$\tilde{E}$};
\draw (228.45,112.2) node [anchor=north west][inner sep=0.75pt]  [color={rgb, 255:red, 0; green, 0; blue, 0 }  ,opacity=1 ]  {$-n$};
\draw (329,138.8) node [anchor=north west][inner sep=0.75pt]    {$\tilde{F}$};
\draw (187,7.3) node [anchor=north west][inner sep=0.75pt]    {$\text{Bl}_{P}(\mathbb{F}_{n})$};
\draw (292.5,35.4) node [anchor=north west][inner sep=0.75pt]    {$-1$};
\draw (286.5,149.4) node [anchor=north west][inner sep=0.75pt]    {$-1$};
\draw (330.5,41.9) node [anchor=north west][inner sep=0.75pt]    {$\hat{E}$};
\draw (307,305.4) node [anchor=north west][inner sep=0.75pt]    {$\alpha _{P}$};
\draw (153.5,165.5) node [anchor=north west][inner sep=0.75pt]    {$\pi $};
\draw (458.5,149.5) node [anchor=north west][inner sep=0.75pt]    {$\pi '$};

\end{tikzpicture}
}
\caption{Elementary transformation \\ $\alpha_P \colon \F_n \dash \F_{n-1}$.}
\label{fig:linkII1.1}
\end{subfigure}
\begin{subfigure}{0.48\textwidth}
\resizebox{1.0\textwidth}{!}{

\tikzset{every picture/.style={line width=0.75pt}} 

\begin{tikzpicture}[x=0.75pt,y=0.75pt,yscale=-1,xscale=1]

\draw   (75.5,264.1) -- (236.3,264.1) -- (236.3,413.1) -- (75.5,413.1) -- cycle ;
\draw [color={rgb, 255:red, 208; green, 2; blue, 27 }  ,draw opacity=1 ][fill={rgb, 255:red, 208; green, 2; blue, 27 }  ,fill opacity=1 ][line width=1.5]    (74.45,338.9) -- (237.35,338.3) ;
\draw [color={rgb, 255:red, 126; green, 211; blue, 33 }  ,draw opacity=1 ][line width=1.5]    (175.4,262.5) -- (175.4,413) ;
\draw  [fill={rgb, 255:red, 144; green, 19; blue, 254 }  ,fill opacity=1 ] (172.2,337.75) .. controls (172.2,335.98) and (173.63,334.55) .. (175.4,334.55) .. controls (177.17,334.55) and (178.6,335.98) .. (178.6,337.75) .. controls (178.6,339.52) and (177.17,340.95) .. (175.4,340.95) .. controls (173.63,340.95) and (172.2,339.52) .. (172.2,337.75) -- cycle ;
\draw   (423.5,264.1) -- (584.3,264.1) -- (584.3,413.1) -- (423.5,413.1) -- cycle ;
\draw [color={rgb, 255:red, 208; green, 2; blue, 27 }  ,draw opacity=1 ][fill={rgb, 255:red, 208; green, 2; blue, 27 }  ,fill opacity=1 ][line width=1.5]    (422.45,338.9) -- (585.35,338.3) ;
\draw [color={rgb, 255:red, 144; green, 19; blue, 254 }  ,draw opacity=1 ][fill={rgb, 255:red, 144; green, 19; blue, 254 }  ,fill opacity=1 ][line width=1.5]    (523.4,262.5) -- (523.4,413) ;
\draw  [fill={rgb, 255:red, 126; green, 211; blue, 33 }  ,fill opacity=1 ] (520.2,299.25) .. controls (520.2,297.48) and (521.63,296.05) .. (523.4,296.05) .. controls (525.17,296.05) and (526.6,297.48) .. (526.6,299.25) .. controls (526.6,301.02) and (525.17,302.45) .. (523.4,302.45) .. controls (521.63,302.45) and (520.2,301.02) .. (520.2,299.25) -- cycle ;
\draw   (243,39.1) -- (403.8,39.1) -- (403.8,188.1) -- (243,188.1) -- cycle ;
\draw [color={rgb, 255:red, 208; green, 2; blue, 27 }  ,draw opacity=1 ][fill={rgb, 255:red, 208; green, 2; blue, 27 }  ,fill opacity=1 ][line width=1.5]    (241.95,113.9) -- (404.85,113.3) ;
\draw [color={rgb, 255:red, 144; green, 19; blue, 254 }  ,draw opacity=1 ][line width=1.5]    (318.7,73.4) .. controls (348.7,118.9) and (358.7,136.4) .. (328.2,171.4) ;
\draw [color={rgb, 255:red, 126; green, 211; blue, 33 }  ,draw opacity=1 ][line width=1.5]    (334.2,46.9) .. controls (353.7,80.4) and (349.2,92.4) .. (320.2,106.4) ;
\draw    (224.2,149.6) -- (160.05,229.94) ;
\draw [shift={(158.8,231.5)}, rotate = 308.61] [color={rgb, 255:red, 0; green, 0; blue, 0 }  ][line width=0.75]    (10.93,-3.29) .. controls (6.95,-1.4) and (3.31,-0.3) .. (0,0) .. controls (3.31,0.3) and (6.95,1.4) .. (10.93,3.29)   ;
\draw    (443,144.1) -- (501.48,210.99) ;
\draw [shift={(502.8,212.5)}, rotate = 228.84] [color={rgb, 255:red, 0; green, 0; blue, 0 }  ][line width=0.75]    (10.93,-3.29) .. controls (6.95,-1.4) and (3.31,-0.3) .. (0,0) .. controls (3.31,0.3) and (6.95,1.4) .. (10.93,3.29)   ;
\draw  [dash pattern={on 4.5pt off 4.5pt}]  (278.8,344.1) -- (398.8,344.1) ;
\draw [shift={(400.8,344.1)}, rotate = 180] [color={rgb, 255:red, 0; green, 0; blue, 0 }  ][line width=0.75]    (10.93,-3.29) .. controls (6.95,-1.4) and (3.31,-0.3) .. (0,0) .. controls (3.31,0.3) and (6.95,1.4) .. (10.93,3.29)   ;

\draw (590.5,332.4) node [anchor=north west][inner sep=0.75pt]  [color={rgb, 255:red, 0; green, 0; blue, 0 }  ,opacity=1 ]  {$E'$};
\draw (427,343.4) node [anchor=north west][inner sep=0.75pt]  [color={rgb, 255:red, 0; green, 0; blue, 0 }  ,opacity=1 ]  {$-( n+1)$};
\draw (506.5,266.4) node [anchor=north west][inner sep=0.75pt]  [color={rgb, 255:red, 0; green, 0; blue, 0 }  ,opacity=1 ]  {$0$};
\draw (516.5,245.4) node [anchor=north west][inner sep=0.75pt]    {$F'$};
\draw (501,289.9) node [anchor=north west][inner sep=0.75pt]    {$Q$};
\draw (409,244.4) node [anchor=north west][inner sep=0.75pt]    {$\mathbb{F}_{n+1}$};
\draw (242,332.4) node [anchor=north west][inner sep=0.75pt]  [color={rgb, 255:red, 0; green, 0; blue, 0 }  ,opacity=1 ]  {$E$};
\draw (79,339.9) node [anchor=north west][inner sep=0.75pt]  [color={rgb, 255:red, 0; green, 0; blue, 0 }  ,opacity=1 ]  {$-n$};
\draw (158.5,266.4) node [anchor=north west][inner sep=0.75pt]  [color={rgb, 255:red, 0; green, 0; blue, 0 }  ,opacity=1 ]  {$0$};
\draw (168,245.4) node [anchor=north west][inner sep=0.75pt]    {$F$};
\draw (156,316.9) node [anchor=north west][inner sep=0.75pt]    {$P$};
\draw (61,243.9) node [anchor=north west][inner sep=0.75pt]    {$\mathbb{F}_{n}$};
\draw (408,98.4) node [anchor=north west][inner sep=0.75pt]  [color={rgb, 255:red, 0; green, 0; blue, 0 }  ,opacity=1 ]  {$\tilde{E}$};
\draw (249.45,120.3) node [anchor=north west][inner sep=0.75pt]  [color={rgb, 255:red, 0; green, 0; blue, 0 }  ,opacity=1 ]  {$-( n+1)$};
\draw (351,47.4) node [anchor=north west][inner sep=0.75pt]    {$\tilde{F}$};
\draw (207,15.4) node [anchor=north west][inner sep=0.75pt]    {$\text{Bl}_{P}(\mathbb{F}_{n})$};
\draw (312.5,43.5) node [anchor=north west][inner sep=0.75pt]    {$-1$};
\draw (306.5,157.5) node [anchor=north west][inner sep=0.75pt]    {$-1$};
\draw (353.5,138.5) node [anchor=north west][inner sep=0.75pt]    {$\hat{E}$};
\draw (327,313.5) node [anchor=north west][inner sep=0.75pt]    {$\alpha _{P}$};
\draw (173.5,173.6) node [anchor=north west][inner sep=0.75pt]    {$\pi $};
\draw (478.5,157.6) node [anchor=north west][inner sep=0.75pt]    {$\pi '$};

\end{tikzpicture}

}
\caption{Elementary transformation \\ $\alpha_P \colon \F_n \dash \F_{n+1}$.}
\label{fig:linkII1.2}
\end{subfigure}

\caption{Sarkisov link of type II.}
\label{fig:linkII}
\end{figure}
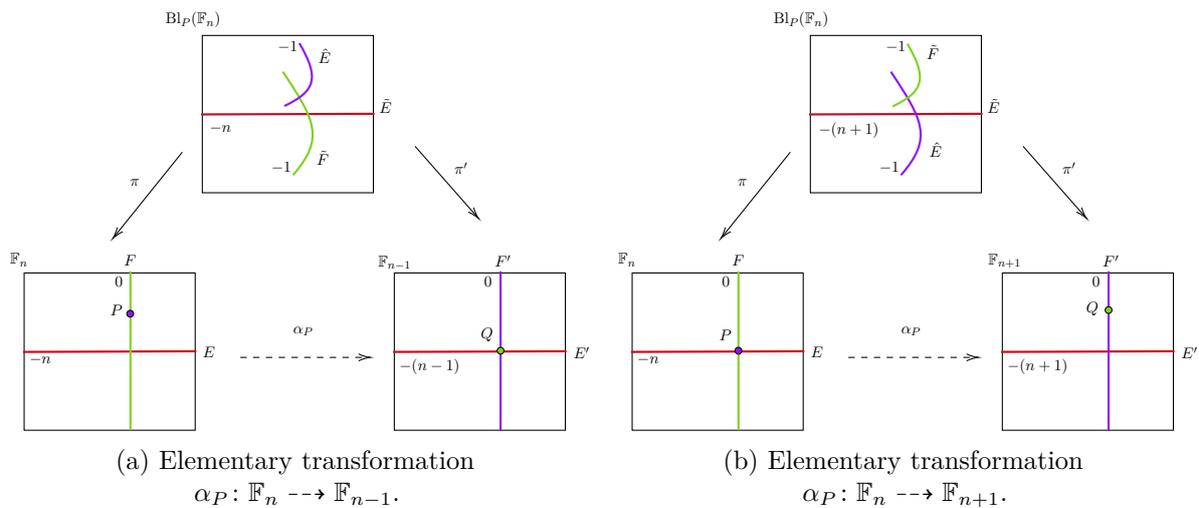

\vspace{1.0cm}

\item A \textit{Sarkisov link of type III} (the inverse of a link of type I) is a commutative diagram:

\begin{center}
\begin{tikzcd}
	{\mathbb{F}_1} \\
	{\mathbb{P}^1} & {\mathbb{P}^2} \\
	& {\text{Spec}(\mathbb{C})}
	\arrow[from=1-1, to=2-1]
	\arrow["\pi", from=1-1, to=2-2]
	\arrow[from=2-2, to=3-2]
	\arrow[from=2-1, to=3-2]
\end{tikzcd}
\end{center}

\noindent where $\pi \colon \F_1 \rightarrow \p^2$ is the blowdown of the negative section of $\F_1$.

\item A \textit{Sarkisov link of type IV} is a commutative diagram:

\begin{center}
\begin{tikzcd}
	{\mathbb{P}^1 \times \mathbb{P}^1} && {\mathbb{P}^1 \times \mathbb{P}^1} \\
	{\mathbb{P}^1} && {\mathbb{P}^1} \\
	& {\text{Spec}(\mathbb{C})}
	\arrow["\tau", from=1-1, to=1-3]
	\arrow["{p_1}"', from=1-1, to=2-1]
	\arrow["{p_2}", from=1-3, to=2-3]
	\arrow[from=2-1, to=3-2]
	\arrow[from=2-3, to=3-2]
\end{tikzcd}
\end{center}

\noindent where $\tau \colon \p^1 \times \p^1 \rightarrow \p^1 \times \p^1$ is the involution which exchanges the two factors.

\end{enumerate}

The proof of the Sarkisov program in the surface case is based on the untwisting of the \textit{Sarkisov degree}. See \cite[Theorem 2.24]{cks} for more details and \cite[Flowchart 1-8-12]{mat} for an explicit flowchart with an equivalent approach. The aim of this data is to measure the complexity of a birational map between the surfaces involved by comparing the linear system associated with the canonical class of the source.

In what follows, let $\bF$ denote either $\p^2$ or some Hirzebruch surface $\F_n$, for $n \geq 0$. 

\begin{defn}\label{defn sark deg dim 2}
    The \textit{Sarkisov degree} of a rational map $\bF \dash \p^2$ given
by the mobile linear system $\Gamma$ is defined as
\begin{enumerate}
    \item $\dfrac{d}{3}$ in case $\bF = \p^2$ and $\Gamma \subset |dH|$, where $H$ is a general line of $\p^2$; or
    \item $\dfrac{b}{2}$ in case $\bF = \F_n $  and $\Gamma \subset |aF+bE|$. (Notice that if $\bF = \p^1 \times \p^1$, the Sarkisov degree is only defined in terms of a choice of one of the two projections $\bF \rightarrow \p^1$.)
\end{enumerate}
\end{defn}

\begin{defn}\label{defn vp sark link}
    A \textit{volume preserving Sarkisov link} is a Sarkisov link previously described with additional data and the following property. There exist divisors $D$ on $\p^2$ and $D'$ on $\F_n$, making $(\p^2,D)$ and $(\F_n,D')$ (t,lc) Calabi-Yau pairs, and all the birational maps that constitute the Sarkisov link are volume preserving for these Calabi-Yau pairs.
\end{defn}

We end this section stating the result of Corti \& Kaloghiros, which holds in all dimensions.

\begin{thm}[cf. \cite{ck} Theorem 1.1]\label{thm vp sark prog}
    Any volume preserving map between Mori fibered Calabi-Yau pairs is a composition of volume preserving Sarkisov links.
\end{thm}
\[\begin{tikzcd}[column sep=scriptsize]
	{(X,D_X)=(X_0,D_0)} & {(X_1,D_1)} & \cdots & {(X_{m-1},D_{m-1})} & {(X_m,D_m)=(Y,D_Y)} \\
	{S=Y_0} & {Y_1} && {Y_{m-1}} & {Y_m=T}
	\arrow["{\phi_1}", dashed, from=1-1, to=1-2]
	\arrow["{\phi_2}", dashed, from=1-2, to=1-3]
	\arrow["{\phi_{m-1}}", dashed, from=1-3, to=1-4]
	\arrow["{\phi_m}", dashed, from=1-4, to=1-5]
	\arrow[from=1-1, to=2-1]
	\arrow[from=1-2, to=2-2]
	\arrow[from=1-4, to=2-4]
	\arrow[from=1-5, to=2-5]
	\arrow["\phi", curve={height=-24pt}, dashed, from=1-1, to=1-5]
\end{tikzcd}\]

Here $\phi$ stands for a volume preserving map between the Mori fibered Calabi-Yau pairs $(X,D_X)/S$ and $(Y,D_Y)/T$, and $\phi_i$ for a volume preserving Sarkisov link in its decomposition.

\section{The decomposition group of a nonsingular plane cubic}\label{dec group cubic}

\paragraph{Decomposition and inertia groups.} These groups were introduced in \cite{giz} in a general context involving the language of schemes. This terminology has its origin in concepts of commutative algebra with some arithmetic implications. Restricting ourselves to the category of projective varieties and rational maps, these groups have the following definitions:

\begin{defn}
    Let $X$ be a projective variety and $\Bir(X)$ be its group of birational automorphisms. Given $Y \subset X$ an (irreducible) subvariety and $\xi_Y$ its generic point, the \textit{decomposition group} of $Y$ in $\Bir(X)$ is the group

    \begin{center}
        $\Bir(X,Y) = \{ \varphi \in \Bir(X)~|~\xi_Y \in \Dom(\varphi) \cap \Dom(\varphi^{-1}),~\varphi(\xi_Y)=\xi_Y \}.$
    \end{center}
\end{defn}

Thus, given $\varphi$ in $\Bir(X,Y)$, one has $\varphi(Y) \subset Y$ and $\varphi |_Y \colon Y \dash Y$ is birational. 

The \textit{inertia group} of $Y$ in $\Bir(X)$ is the group 

\begin{center}
$\{ \varphi \in \Bir(X,Y)~|~\varphi |_Y = \Id_Y~\text{as birational map} \}$.
\end{center}

When $X$ is normal and $Y$ is a prime divisor $D$ such that $(X,D)$ is a Calabi-Yau pair, one has $\Bir^{\vp}(X,D)=\Bir(X,D)$ provided the pair $(X,D)$ has canonical singularities. This is the content of the following result: 

\begin{prop}[cf. \cite{acm} Proposition 2.6]\label{vp can pair} Let $(X,D_X)$ and $(Y,D_Y)$ be (t,c) Calabi-Yau pairs, and
  $f\colon X\dasharrow Y$ an arbitrary birational map.  Then
  $f\colon (X,D_X)\dasharrow(Y,D_Y)$ is volume preserving if and only
  if $f_{*}D_X = D_Y$ and $f^{-1}_{*}D_Y = D_X$. (This condition is
  equivalent to asking that the restriction of $f$ to each component
  of $D_X$ is a birational map to a component of $D_Y$, and the same
  for $f^{-1}$).

In particular, if $(X,D)$ is a (t,c)  Calabi-Yau pair with $D$ irreducible, then $\Bir^{\vp}(X,D)$ coincides with the group of birational self-maps $f\colon X \dasharrow X$ that restrict to birational self-maps $f|_D \colon D \dash D$.
\end{prop}

When the ambient variety $X$ is $\p^n$, we denote such groups by $\Dec(Y)$ and $\Ine(Y)$, respectively.

Let $C \subset \p^2$ be an irreducible nonsingular cubic. We have readily that $(\p^2,C)$ is a Calabi-Yau pair with (t,c) singularities according to the Definition \ref{defn pairs}. In particular, one can easily check that $C$ is an elliptic curve. The following theorem by Pan \cite{pan} gives an interesting property of the elements of $\Dec(C)$:

\vspace{0.01cm}
\begin{thm}[cf. \cite{pan} Théorème 1.3]\label{pan thm 1.3}
Let $C \subset \p^{2}$ be an irreducible, nonsingular and nonrational curve and suppose there exists $\phi \in \Dec(C)\setminus \PGL(3,\C)$. Then $\deg(C)=3$ and $\Bs(\phi) \subset C$, where $\Bs(\phi)$ denotes the set of proper base points of $\phi$.
\begin{figure}[htb]
\centering
\resizebox{0.8\textwidth}{!}%
{\tikzset{every picture/.style={line width=0.75pt}} 
\begin{tikzpicture}[x=0.75pt,y=0.75pt,yscale=-1,xscale=1]
\draw    (313.53,2.87) ;
\draw   (117,34) -- (265.26,34) -- (201.72,164) -- (53.46,164) -- cycle ;
\draw [line width=1.5]    (94.8,128) .. controls (120.8,209) and (209.8,-34.29) .. (216,80.71) ;
\draw  [fill={rgb, 255:red, 144; green, 19; blue, 254 }  ,fill opacity=1 ] (120,143) .. controls (120,140.79) and (118.21,139) .. (116,139) .. controls (113.79,139) and (112,140.79) .. (112,143) .. controls (112,145.21) and (113.79,147) .. (116,147) .. controls (118.21,147) and (120,145.21) .. (120,143) -- cycle ;
\draw  [fill={rgb, 255:red, 144; green, 19; blue, 254 }  ,fill opacity=1 ] (135,123.71) .. controls (135,121.51) and (136.79,119.71) .. (139,119.71) .. controls (141.21,119.71) and (143,121.51) .. (143,123.71) .. controls (143,125.92) and (141.21,127.71) .. (139,127.71) .. controls (136.79,127.71) and (135,125.92) .. (135,123.71) -- cycle ;
\draw  [fill={rgb, 255:red, 144; green, 19; blue, 254 }  ,fill opacity=1 ] (152.5,101.4) .. controls (152.5,99.38) and (154.14,97.74) .. (156.16,97.74) .. controls (158.18,97.74) and (159.82,99.38) .. (159.82,101.4) .. controls (159.82,103.42) and (158.18,105.06) .. (156.16,105.06) .. controls (154.14,105.06) and (152.5,103.42) .. (152.5,101.4) -- cycle ;
\draw  [fill={rgb, 255:red, 144; green, 19; blue, 254 }  ,fill opacity=1 ] (165.63,81.2) .. controls (165.63,79.18) and (167.27,77.54) .. (169.29,77.54) .. controls (171.31,77.54) and (172.95,79.18) .. (172.95,81.2) .. controls (172.95,83.22) and (171.31,84.86) .. (169.29,84.86) .. controls (167.27,84.86) and (165.63,83.22) .. (165.63,81.2) -- cycle ;
\draw  [fill={rgb, 255:red, 144; green, 19; blue, 254 }  ,fill opacity=1 ] (198.9,49.67) .. controls (198.9,47.65) and (200.54,46.01) .. (202.56,46.01) .. controls (204.58,46.01) and (206.22,47.65) .. (206.22,49.67) .. controls (206.22,51.69) and (204.58,53.33) .. (202.56,53.33) .. controls (200.54,53.33) and (198.9,51.69) .. (198.9,49.67) -- cycle ;
\draw  [dash pattern={on 0.84pt off 2.51pt}]  (113.33,144.67) -- (87.33,146.4) ;
\draw  [fill={rgb, 255:red, 126; green, 211; blue, 33 }  ,fill opacity=1 ] (83.67,146.4) .. controls (83.67,144.37) and (85.31,142.73) .. (87.33,142.73) .. controls (89.36,142.73) and (91,144.37) .. (91,146.4) .. controls (91,148.43) and (89.36,150.07) .. (87.33,150.07) .. controls (85.31,150.07) and (83.67,148.43) .. (83.67,146.4) -- cycle ;
\draw  [dash pattern={on 0.84pt off 2.51pt}]  (172.33,77.67) -- (190,51.67) ;
\draw  [fill={rgb, 255:red, 126; green, 211; blue, 33 }  ,fill opacity=1 ] (186.8,51.32) .. controls (186.8,49.37) and (188.39,47.78) .. (190.34,47.78) .. controls (192.3,47.78) and (193.88,49.37) .. (193.88,51.32) .. controls (193.88,53.28) and (192.3,54.87) .. (190.34,54.87) .. controls (188.39,54.87) and (186.8,53.28) .. (186.8,51.32) -- cycle ;
\draw   (431.5,34.5) -- (579.76,34.5) -- (516.22,164.5) -- (367.96,164.5) -- cycle ;
\draw [line width=1.5]    (410.8,127) .. controls (436.8,208) and (525.8,-35.29) .. (532,79.71) ;
\draw  [dash pattern={on 4.5pt off 4.5pt}]  (261.7,104.7) -- (362.2,104.7) ;
\draw [shift={(364.2,104.7)}, rotate = 180] [color={rgb, 255:red, 0; green, 0; blue, 0 }  ][line width=0.75]    (10.93,-3.29) .. controls (6.95,-1.4) and (3.31,-0.3) .. (0,0) .. controls (3.31,0.3) and (6.95,1.4) .. (10.93,3.29)   ;
\draw  [dash pattern={on 0.84pt off 2.51pt}]  (86.32,125.88) -- (87.33,142.73) ;
\draw  [fill={rgb, 255:red, 245; green, 166; blue, 35 }  ,fill opacity=1 ] (82.86,122.42) .. controls (82.86,120.51) and (84.41,118.96) .. (86.32,118.96) .. controls (88.23,118.96) and (89.78,120.51) .. (89.78,122.42) .. controls (89.78,124.33) and (88.23,125.88) .. (86.32,125.88) .. controls (84.41,125.88) and (82.86,124.33) .. (82.86,122.42) -- cycle ;
\draw  [fill={rgb, 255:red, 144; green, 19; blue, 254 }  ,fill opacity=1 ] (37,197.3) .. controls (37,195.2) and (38.7,193.5) .. (40.8,193.5) .. controls (42.9,193.5) and (44.6,195.2) .. (44.6,197.3) .. controls (44.6,199.4) and (42.9,201.1) .. (40.8,201.1) .. controls (38.7,201.1) and (37,199.4) .. (37,197.3) -- cycle ;
\draw  [fill={rgb, 255:red, 126; green, 211; blue, 33 }  ,fill opacity=1 ] (21.33,220.85) .. controls (21.33,218.94) and (22.88,217.4) .. (24.78,217.4) .. controls (26.69,217.4) and (28.23,218.94) .. (28.23,220.85) .. controls (28.23,222.76) and (26.69,224.3) .. (24.78,224.3) .. controls (22.88,224.3) and (21.33,222.76) .. (21.33,220.85) -- cycle ;
\draw  [fill={rgb, 255:red, 245; green, 166; blue, 35 }  ,fill opacity=1 ] (39.36,220.92) .. controls (39.36,219.01) and (40.91,217.46) .. (42.82,217.46) .. controls (44.73,217.46) and (46.28,219.01) .. (46.28,220.92) .. controls (46.28,222.83) and (44.73,224.38) .. (42.82,224.38) .. controls (40.91,224.38) and (39.36,222.83) .. (39.36,220.92) -- cycle ;

\draw (200,84.9) node [anchor=north west][inner sep=0.75pt]    {$C$};
\draw (261,9.9) node [anchor=north west][inner sep=0.75pt]    {$\mathbb{P}^{2}$};
\draw (516,83.9) node [anchor=north west][inner sep=0.75pt]    {$C$};
\draw (577,8.9) node [anchor=north west][inner sep=0.75pt]    {$\mathbb{P}^{2}$};
\draw (306,74.4) node [anchor=north west][inner sep=0.75pt]    {$\phi $};
\draw (57,188) node [anchor=north west][inner sep=0.75pt]   [align=left] {{\fontfamily{helvet}\selectfont proper base points of} $\displaystyle \phi $};
\draw (56,210) node [anchor=north west][inner sep=0.75pt]   [align=left] {possible infinitely near base points of $\displaystyle \phi $};

\end{tikzpicture}
    }  

\caption{Element of $\Dec(C)$.}
\label{fig:dec}
\end{figure}
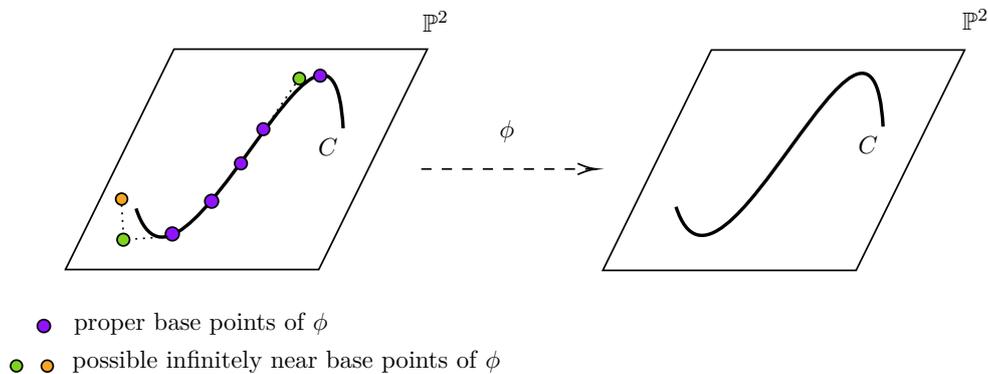
\end{thm}

By Proposition \ref{vp can pair}, $\Bir^{\vp}(\p^2,C) = \Dec (C)$, that is, the elements of $\Dec(C)$ are exactly the volume preserving self-maps of the Calabi-Yau pair $(\p^2,C)$. \ref{vp can pair}

Under the point of view of log Calabi-Yau geometry, Theorem \ref{thm vp sark prog} ensures the existence of a volume preserving factorization for any element of $\Dec(C)$. On the other hand, the elements of $\Dec(C)$ can be seen as ordinary maps in $\Bir(\p^2)$ and consequently they admit a Sarkisov factorization. There is no reason at first for this standard Sarkisov factorization to be volume preserving. By the adjective \textit{standard} here, we mean a Sarkisov factorization obtained by running the usual Sarkisov program without taking into account the volume preserving property. Our result Theorem \ref{thm sark=vp cubic} says that the Sarkisov algorithm in dimension 2 is automatically volume preserving for an element of $\Dec(C)$.

\begin{defn}
Given $\Gamma$ a linear system on a nonsingular projective surface $S$, the \textit{multiplicity} $m_P$ of $\Gamma$ at a point $P \in S$ is the multiplicity of a general member there, that is,
\begin{center}
    $m_P \coloneqq \min\{\mult_P(C);\ C \in \Gamma\}$.
\end{center}
\end{defn}

\begin{nota} We denote by $\Bs(\phi)$ the proper base locus of a rational map $\phi \colon X \dash Y$ between projective varieties, and by $\uBs(\phi)$ its full base locus, including the infinitely near one. See \cite[Section 2.2]{cks} for more details about this notion in the surface case. Moreover, there exists a natural partial ordering of the points proper or infinitely near to any nonsingular projective surface. We write $P \prec Q$ if $Q$ is infinitely near to $P$.
\end{nota}

Before stating and proving Theorem \ref{thm sark=vp cubic}, we will show a couple of lemmas followed by a stronger version of Theorem \ref{pan thm 1.3}.

\begin{lem}\label{blowup vp}
    Let $S$ be a nonsingular projective surface admitting a Calabi-Yau pair structure $(S,C)$ with boundary divisor $C$ nonsingular. Consider $P \in S$. Then $f \colon (\Bl_P(S),\Tilde{C}) \rightarrow (S,C)$ is volume preserving if and only if $P \in C$.
\end{lem}

\begin{proof}
    Let $E \coloneqq \Exc(f)$. By the adjunction formula, we can write $K_{\Bl_P(S)} = f^*K_S + E$. Since $C$ is nonsingular, we have $\Tilde{C} = f^*C -mE$, where $m \in \{0,1\}$. Summing up we get 
    \begin{center}
        $K_{\Bl_P(S)} + \Tilde{C} = f^*(K_S +C) +(1-m)E$.
    \end{center}
    
    Since $f_*\Tilde{C}=C$, by Proposition \ref{prop crep bir morp}, $(\Bl_P(S),\Tilde{C})$ is a Calabi-Yau pair and $f$ is volume preserving if and only if $1-m = 0$, that is, if and only if $P \in C$.
\end{proof}

\begin{lem}\label{blowndown vp}
     Let $S$ be a nonsingular projective surface admitting a Calabi-Yau pair structure $(S,C)$ with boundary divisor $C$ nonsingular. Let $E \subset S$ be a $(-1)$-curve, $f \colon S \rightarrow S'$ the contraction of $E$ and $C'=f(C) \subset S'$. Then  $f \colon (S,C) \rightarrow (S',C')$ is volume preserving if and only if $C \cdot E = 1$ (and therefore $C'$ is nonsingular).
\end{lem}

\begin{proof}
By the Castelnuovo contratibility theorem, we have $E = \Exc(f)$ and $S \simeq \Bl_{P}(S')$ for $\{P\} = f(E)$. Since $C \neq E$, we have that $C'$ is indeed a curve and $C = \Tilde{C'}$.

Set $m \coloneqq m_P(C')$. Observe that $C = f^*C' - mE$ and so $C \cdot E = m$.

Notice that
\begin{center}
    $0 = f_*(K_S + C) = K_{S'}+C'$,
\end{center}
and therefore $(S',C')$ is a Calabi-Yau pair.

By analogous computations to the previous lemma, we get
\begin{center}
    $K_S+C = f^*(K_{S'}+C') + (1-m)E$.
\end{center}

Since $f_*C=C'$, by Proposition \ref{prop crep bir morp}, $f$ is volume preserving if and only if $1-m = 0$, that is, if and only if $P \in C'$. Furthermore, this also shows that $C'$ is nonsingular at $P$ and therefore $C'$ is nonsingular. 
\end{proof}

The following lemma says that not all Mori fibered Calabi-Yau pairs in dimension 2 admit an irreducible boundary divisor:

\begin{lem}\label{lem Mf CY pairs in dim 2}
The only (rational) Mori fibered spaces in dimension 2 that admit an irreducible divisor $D$ such that $(S,D)$ is a Calabi-Yau pair are:
\begin{center}
    $\p^2/\Spec(\C)$ and $\F_n/\p^1$ for $n \in \{0,1,2\}$.
\end{center}

Moreover, if $D$ is nonsingular, then $D$ is isomorphic to an elliptic curve.
\end{lem} 

\begin{proof} The case $S=\p^2$ is immediate. For the remaining cases $\F_n$, this result follows from intersecting $D$ with the negative section.

If $(\F_n,D)$ is a Calabi-Yau pair, then $D \sim (n+2)F+2E$ because $K_{\F_n}=-(n+2)F-2E$.

Since $D \neq E$, we must then have that $D \cdot E \geq 0$. Observe that
\begin{center}
    $D \cdot E = ((n+2)F+2E) \cdot E = (n+2)-2n=-n+2 \geq 0 \Leftrightarrow n \leq 2$.
\end{center}

We have just shown that if $(\F_n,D)$ is a Calabi-Yau pair, then $n\leq 2$. For the converse, it is easy to construct examples of such pairs by applying conveniently the Lemmas \ref{blowup vp} \& \ref{blowndown vp} to the Calabi-Yau pair $(\p^2,D)$, where $D$ is a nonsingular cubic.

Besides being irreducible, assume $D$ is nonsingular. By the adjunction formula, we have that $g=1$. Therefore, $D$ is an elliptic curve and the result then follows.
\end{proof}
Let $C \subset \p^2$ be a nonsingular cubic. Recall that Pan \cite{pan} showed that if $\phi \in \Dec(C) \setminus \PGL(3,\C)$, then $\Bs(\phi) \subset C$. Using the volume preserving version of the Sarkisov program, it is possible to show that more is true: $\uBs(\phi)$ is contained in $C$. Of course, this is an abuse of language since the points in $\uBs(\phi)$ do not lie on $\p^2$, but on some infinitesimal neighborhood of the points in $\Bs(\phi)$.

Thus, $\uBs(\phi)$ is contained in $C$ means that the points in $\uBs(\phi)$ belong to the strict transforms of $C$ intersected with the infinitesimal neighborhoods of the points in $\Bs(\phi)$. The subsequent lemma presents a stronger form of Theorem \ref{pan thm 1.3}.

\begin{lem}\label{lem inf base locus cont in C}
    Let $C \subset \p^2$ be a nonsingular cubic. Consider $\phi \in \Dec(C) \setminus \PGL(3,\C)$. Then $\uBs(\phi) \subset C$.
\end{lem}

\begin{proof}
    Let $P_1 \prec \cdots \prec P_l$ be any maximal increasing sequence of base points in $\uBs(\phi)$, starting with a base point $P_1 \in \Bs(\phi)$. We will show by increasing induction on $i$ that all $P_i \in C$. The base case $i=1$ follows from Theorem \ref{pan thm 1.3}, which shows that $\Bs(\phi) \subset C$.

    Consider a volume preserving Sarkisov factorization of $\phi$:
\[\begin{tikzcd}
	{(\mathbb{P}^2,C)=(S_0,C_0)} & {(S_1,C_1)} & \cdots & {(S_k,C_k)=(\mathbb{P}^2,C)}
	\arrow["{\phi_0}", dashed, from=1-1, to=1-2]
	\arrow["{\phi_{k-1}}", dashed, from=1-3, to=1-4]
	\arrow["{\phi_1}", dashed, from=1-2, to=1-3]
	\arrow["\phi", curve={height=-18pt}, dashed, from=1-1, to=1-4]
\end{tikzcd}.\]

Observe that each $S_j \simeq \p^2$ or $\F_n$ for $n \in \{0,1,2\}$ by Lemma \ref{lem Mf CY pairs in dim 2}. Indeed, \cite[Lemma 2.8]{acm} establishes that canonicity is retained for volume preserving maps between Mori fibered Calabi-Yau pairs. Since $(\p^2,C)$ is a (t,c) Calabi-Yau pair, this implies that all the intermediate ones $(S_j,C_j)$ are also (t,c). By Proposition \ref{vp can pair}, we have that each $C_j$ is the strict transform of $C$ on $S_j$. Since each $(S_j,C_j)$ is, in particular, a dlt Calabi-Yau pair, by \cite[Proposition 5.51]{km} or \cite[Remark 1.3(1)]{acm}, we have that each $C_i$ is normal, therefore nonsingular. Thus, it follows that $C_j \simeq C$ for all $j \in \{ 0,1,\ldots,k-1 \}$.

Notice that a point in $\uBs(\phi)$ appears in the proper base locus of some induced birational map $\phi \circ \phi_0^{-1} \circ \ldots \phi_j^{-1}$ when one runs the Sarkisov algorithm. Thus, the key observation is that for all $i \in \{0,1,\ldots,l \}$, there exists $j_i$ such that a Sarkisov link of type I or II $\phi_{j_i} \colon S_{j_i} \dash S_{j_i+1}$ starts with the blowup of the image of $P_i$ on $S_{j_i}$. 

Since the Sarkisov link is volume preserving, in particular so is the blowup of $P_i$ that initiates the link by Definition \ref{defn vp sark link}. Conclude from Lemma \ref{blowup vp}. 
  
\end{proof}

\begin{rmk}\label{rmk on the lem inf base locus cont in C}
    The previous Lemma \ref{lem inf base locus cont in C} can also be proved without the use of the volume preserving Sarkisov program. Indeed, one can argue as follows. Consider a (minimal) decomposition of a given birational map in $\Dec(C)$ into blowups and blowdowns and keep track of the boundary throughout the process. Since in the end we have a boundary and not a sub-boundary, blowing up a point that does not belong to the strict transform of $C$ at any intermediate step produces a component in the crepant pullback with negative coefficient. Such a component must necessarily be contracted at a later stage. Hence, we may skip this step of the blowup by using the universal property and obtain a contradiction to the minimality of the decomposition. Note that this argument does not use the assumption on the singularities of the curve. We nevertheless keep the proof of Lemma \ref{lem inf base locus cont in C} as it is in order to illustrate an additional application of the volume preserving Sarkisov program.
\end{rmk}
 
\begin{thm}\label{thm sark=vp cubic}
Let $C \subset \p^2$ be a nonsingular cubic. Any standard Sarkisov program applied to an element of $\Dec(C)$ is automatically volume preserving.
\end{thm}

\begin{proof} 
By Lemma \ref{lem inf base locus cont in C} and Remark \ref{rmk on the lem inf base locus cont in C}, it suffices to check that a Sarkisov link which starts with a blowup (or the inverse of a blowup) has the blownup
point lying on the boundary divisor, which follows from the arguments therein. The remainder follows from the fact that the pushforward of a Calabi-Yau pair under a contraction is again a Calabi-Yau pair, and a nonsingular boundary curve cannot be mapped to a singular one in this setting, as explained in Lemma \ref{blowndown vp}. Combining all of this together, we get the result.

\end{proof}

\begin{rmk}
The canonicity of the Calabi-Yau pairs in the previous results was really necessary. Consider the divisor $D = L_1 + L_2 + L_3$ on $\p^2$ given by the sum of the three coordinate lines. The Calabi-Yau pair $(\p^2,D)$ is strictly log canonical, that is, it is (t,lc) but not (t,c). \footnote{Usually, by strictly log canonical, it is meant log canonical but not Kawamata log terminal (klt). As pointed out right after Definition 3.1, discrepancies with respect to Calabi-Yau pairs are integer numbers. So, klt means canonical under this assumption.} One may check that the standard quadratic transformation $(x : y : z) \mapsto (yz : xz : xy)$ is a volume preserving self-map for this pair. The Calabi-Yau pairs appearing in its volume preserving (standard) Sarkisov decomposition are such that the boundaries are not given by the strict transform of the initial one $L_1+L_2+L_3$.   

The notion of volume preserving birational map is defined not with respect to strict transform of the boundary, but with respect to crepant equivalent pair which is pullback to the resolution and the pushforward of the initial boundary (as in Definition \ref{defn crep bir map}). The point is that when our initial Calabi-Yau pair is (t,c), the Calabi-Yau pairs appearing in a factorization of a volume preserving self-map have the form $(X,D_X)$, where $D_X$ is the strict transform of the initial boundary divisor. This is a consequence of Proposition \ref{vp can pair}.
          
\end{rmk}

The Sarkisov program in dimension 2 can also yield a factorization of elements in $\Bir(\p^2)$ into de Jonquières maps, incorporating additional steps. See \cite[Theorem 2.30]{cks}. Recall that such maps are elements of $\Bir(\p^2)$ that preserve a pencil of lines. In other words, $J \in \Bir(\p^2)$ is a \textit{de Jonquières map} if there exist $P,Q \in \p^2$ such that $J$ sends all the lines through $P$ to lines through $Q$, up to a finite number. We will call such $P,Q$ the centers of de Jonquières map $J$.

\newpage

\begin{figure}[h]
\centering
\includegraphics[scale=0.20]{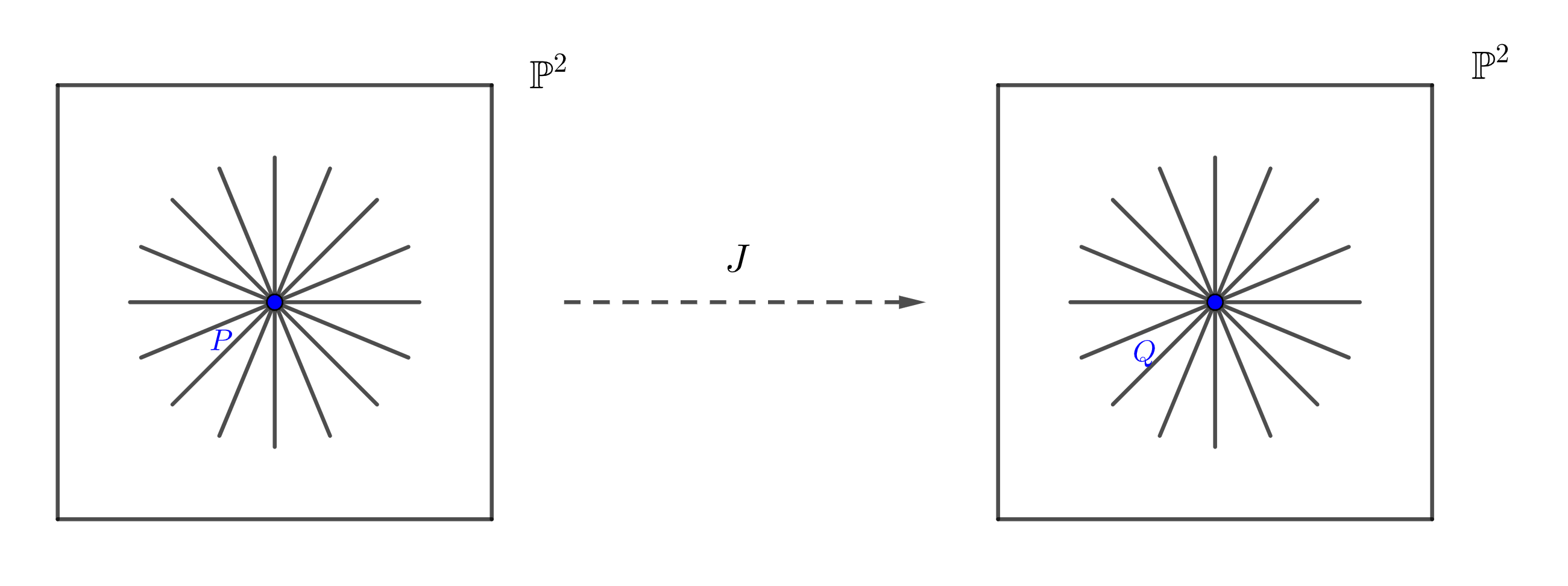}
\caption{de Jonquières map.}
\label{fig:jonqmap}
\end{figure}

\begin{ex}
    The standard quadratic transformation $(x : y : z) \mapsto (yz : xz : xy)$ is a de Jonquières map. Indeed, take $P=Q=(1:0:0)$. One can show the image of any line through $P$ distinct from $\{y =0 \} $ and $\{z =0 \} $ is another line through $P$. 
\end{ex}

\begin{figure}[h]
\centering
\includegraphics[scale=0.20]{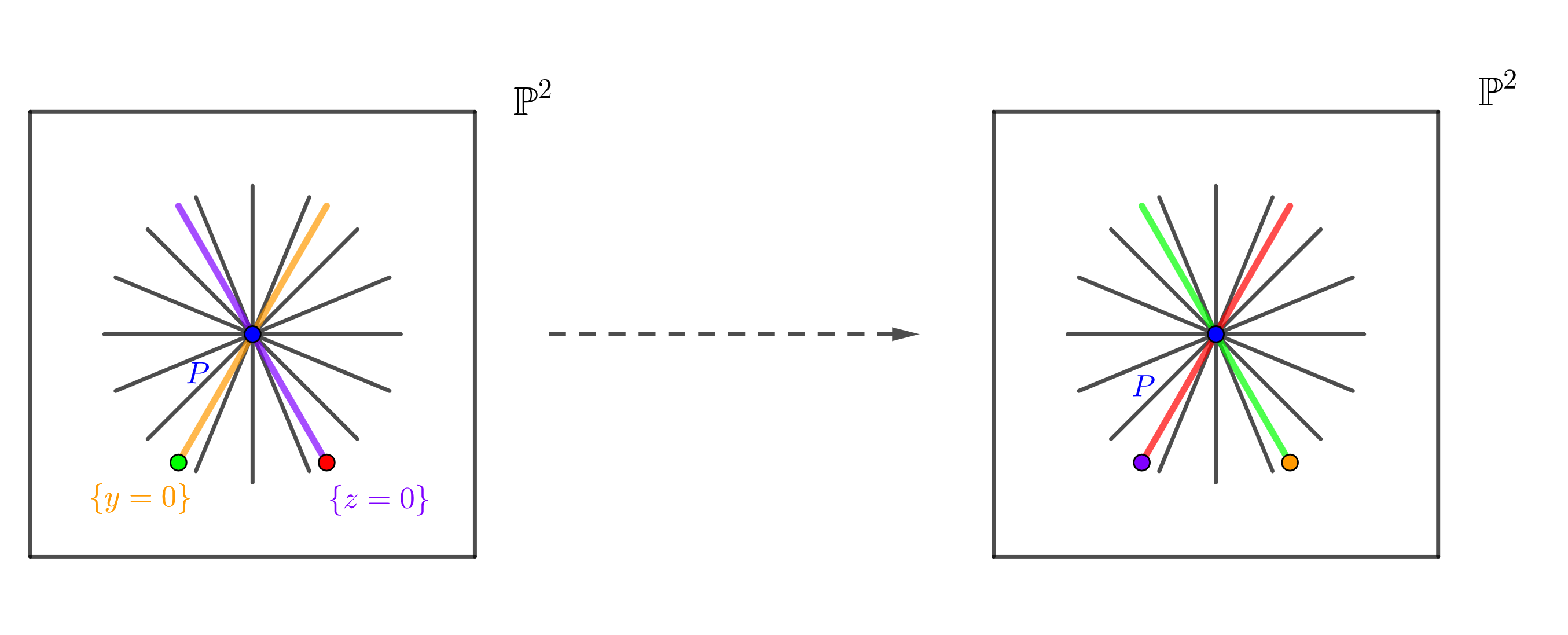}
\caption{Standard quadratic transformation.}
\label{fig:standard}
\end{figure}

The arguments in the proof of \cite[Theorem 2.30]{cks} allow to have the following immediate corollary:

\begin{cor}\label{cor jonq}
    The centers of de Jonquières transformations obtained from the (volume preserving) Sarkisov program applied to any element of $\Dec(C)$ belong to the cubic and its strict transform.
\end{cor}

\begin{proof}
    Given a Sarkisov factorization of an element of $\Bir(\p^2)$, by \cite[Proposition 2.31]{cks}, a de Jonquières one is obtained by possibly incorporating some extra Sarkisov links of types II and III. We are then led to the same setting as in the proof of Theorem \ref{thm sark=vp cubic} and by the same arguments therein, we can assure that the centers of the Jonquières maps that will appear indeed belong to the cubic and its strict transform.
\end{proof}

We recall that Lemma \ref{lem Mf CY pairs in dim 2} restricts possibilities for the Mori fibered spaces appearing in a volume preserving factorization of $\phi \in \Dec(C)$. This is illustrated by examples in \cite[Subsection 4.3.1]{alv}.

\subsection{The canonical complex of $C \subset \p^2$}

Let $C \subset \p^2$ be an irreducible plane curve (not necessarily nonsingular). We have a complex (not necessarily exact) induced by the natural action $\rho$ of $\Dec(C)$ on $C$ 
\begin{equation}\label{cc}
    1 \longrightarrow \Ine(C) \longrightarrow \Dec(C) \overset{\rho} \longrightarrow \Bir(C) \longrightarrow 1,
\end{equation}
where $\Ine(C)$ is identified with $\ker(\rho)$. This complex is called the \textit{canonical complex} of the pair $(\p^2,C)$ and the obstruction to its exactness is the surjectivity of $\rho$. 

In \cite{bpv1}, Blanc, Pan \& Vust studied canonical complexes under the usual trichotomy: genera $g \ge 2, g=1$ and $g=0$. See \cite{bpv1} for interesting examples in which the map $\Dec(C) \rightarrow \Bir(C)$ is not surjective. In the case where $C$ is nonsingular, one has $\Bir(C) = \Aut(C)$. 

If $g \ge 2$, one can easily check that $C$ has degree $>3$. See \cite[Proposition 5, Chapter 8]{ful}. By the first part in the proof of the Theorem \ref{pan thm 1.3} or \cite[Corollaire 3.6]{pan}, in this case, the group $\Dec(C)$ is trivial in the sense that it is given by the automorphisms of $\p^2$ that preserve $C$. That is, if $g \geq 2$, then 
\begin{center}
    $\Dec(C) = \Aut(\p^2,C) \coloneqq \Dec(C) \cap \Aut(\p^2)$.
\end{center}

This fact together with the following result due to Matsumura \& Monsky \cite{mm} when $n \ge 2$ and Chang \cite{ch} when $n=1$ implies that the canonical complex of the pair $(\p^2,C)$ is exact.

\begin{thm}[cf. \cite{mm} Theorem 2 and \cite{ch} Theorem 1] Let $n$ and $d$ be positive integers and $X$ be a nonsingular hypersurface of degree $d$ in $\p^{n+1}$. If $(n,d) \neq (2,4),(1,3)$, then the natural group homomorphism $\Aut(\p^n,X) \rightarrow \Aut(X)$ is surjective.
    
\end{thm}

In the case of $d > n$, this result can also be obtained through adjoint systems \cite[Remarque 3.7]{pan}. When $(n,d)=(1,3)$, then $X = C$ is a nonsingular plane cubic. The result by Pan about its decomposition group indicates the existence of nonlinear maps inducing automorphisms of $C$ by restriction. So $\Dec(C)$ is larger than $\Aut(\p^2,C)$.

In \cite{bl1}, Blanc showed that the inertia group of $C$ is generated by its elements of degree 3, and except for the identity, such elements are the ones with lowest degree. By \cite[Theorem 6]{giz} and \cite[Theorem 2.2]{og} we have that $\Aut(C)$ is no longer derived from $\Aut(\p^2,C)$, but from $\Dec(C)$. Consequently, this implies the exactness of the canonical complex of $(\p^2,C)$. From this, we can therefore also describe $\Aut(C)$ as the quotient group $\Dec(C)/\Ine(C)$.

\paragraph{Expectation by Blanc, Pan \& Vust.} Once we know that a given short complex is exact, we may ask if it splits or not. In \cite{bpv1}, Blanc, Pan \& Vust remarked on the problem of splitting or non-splitting of the canonical complex of $(\p^2,C)$. In this subsection, we confirm their expectation for the latter.

The curve $C$ is elliptic, and therefore it is also an algebraic group. One has $\Aut(C) \simeq C \rtimes \mathbb{Z}_d$, where $C$ is identified with its group of translations and $d \in \{2,4,6\}$, depending on the $j$ invariant of $C$. More precisely, $\mathbb{Z}_d \simeq \Aut(C,O)$, the group of automorphisms of $C$ which fix the neutral element $O$ of the group operation of $C$.

Denote by $\oplus$ the group law on $C$ and fix $O \in C$ a neutral element. Given $P \in C \setminus \{O\}$, consider a map $\varphi_P \in \Dec(C)$ which by restriction to $C$ induces the translation $T_P$ by $P$, that is, $(\varphi_P)|_C(Q)=Q \oplus P$ for all $Q \in C$.

We observe that there are infinitely many maps in $\Dec(C)$ which will induce the same translation $T_P$ on $C$. Any composition with an element of $\Ine(C)$ plays the same role. We point out that $\Ine(C)$ as well as $\Dec(C)$ are infinite uncountable groups. This is based on the existence of a free subgroup of the former, their descriptions in terms of presentations and the cardinality of our ground field $\C$. See \cite[Theorem 6]{bl1} and \cite[Théorèm 1.4]{pan}.  

One can explicitly obtain such a map $\varphi_P$, for example, by homogenizing the expressions of the group law on $C$ with affine coordinates in its Weierstrass normal form, and extending them to $\p^2$ as a rational map. In what follows we make this explicit.

After a suitable change of coordinates, we can write the equation of $C$ in the Weierstrass normal form as 
\begin{center}
    $y^2=x^3+px+q,$
\end{center}

\noindent where $(p,q) \in \C^2 \setminus \{(0,0)\}$, and $(x,y)$ are affine coordinates of $\mathbb{A}^2_{(x,y)} = \{z \neq 0\} \subset \p^2$. The neutral element $O$ of $\oplus$ is the unique intersection of $C$ with the line at infinity: $O=(0:1:0)$, which is an inflection point.

Let us consider now the translation in $C$ by $(a,b) \in C$. In these coordinates, the group law can be explicitly described as $(x',y') = (x,y) \oplus (a,b)$ if and only if
\begin{center}
    $x' = \left( \dfrac{y-b}{x-a} \right)^2 - x-a,~y'=\dfrac{y-b}{x-a}(x'-a)+b$.
\end{center}

Let $\p^2$ have homogeneous coordinates $(x:y:z)$. Consider the rational map $\varphi_P \colon \p^2 \dash \p^2$ defined by the same equations above in the open subset $\{ z \neq 0 \}$. One can verify birationality and that indeed $\varphi \in \Dec(C)$. In this expression, $\varphi_P$ is not defined when $x=a$.

Let us extend $\varphi_P$ to the largest possible open subset of $\p^2$. Performing some algebraic manipulations and homogenizing, the extension obtained, also denoted by $\varphi_P$, becomes
\begin{align*}
\varphi_P \colon \p^2 & \dash \p^2 \\
(x:y:z)&\longmapsto (F_1(x,y,z):F_2(x,y,z):F_3(x,y,z)),
\end{align*}
with 
\begin{align*}
    F_1(x,y,z) & = z(y-bz)^2(x-az)-(x^2-a^2z^2)(x-az)^2 ,\\
    F_2(x,y,z) & = z(y-bz)^3-(y-bz)(x+2az)(x-az)^2,\\
    F_3(x,y,z) & = z(x-az)^3
\end{align*}
for all $(x:y:z) \in \Dom(\varphi_P) \coloneqq \p^2 \setminus \Bs(\varphi_P)$. One can check that $\Bs(\varphi_P) = V(F_1,F_2,F_3) = \{P, O\} \subset C$, as predicted by Theorem \ref{pan thm 1.3}.

Let $\Gamma \subset |4H|$ be the linear system associated to $\varphi_P$, where $H$ denotes a general line of $\p^2$. We have that $\Gamma$ is contained in the linear system of plane quartics passing through $P$ and $O$ with certain multiplicities $m_P$ and $m_O$, respectively. To compute them, let us make use of the Noether-Fano equations or equations of condition from \cite[Section 2.5]{alb} or \cite[Theorem 2.9]{cks}.

These equations imply that we have two possibilities for the multiplicities of base points (including infinitely near ones) of a plane birational map of degree $4$ in nonincreasing order, namely, $(3,1,1,1,1,1,1)$ or $(2,2,2,1,1,1)$.

A careful analysis of the general member of $\Gamma$ shows that we are in the first case, with $m_P=3$ and $m_O=1$. More precisely, $\Gamma$ is contained in the linear system of plane quartics passing through $P$ with multiplicity $3$ and passing through $O$ with multiplicity $1$ and sharing one tangent direction at $O$ and higher order Taylor terms up to order 5. The shared tangent direction is exactly the tangent direction to $C$ at $O$.

Its homaloidal type is $(4;3,1^6)$, where the coordinate before the semicolon indicates its degree and the further ones in nonincreasing order represent the multiplicities of all base points including the infinitely near ones.

According to \cite[Definition 2.28(2)]{cks}, such homaloidal type makes $\varphi_P$ a de Jonquières map and the configuration of the seven base points of $\varphi_P$ (including the infinitely near ones) is as follows
\begin{center}
$P,$ \\
$O \prec P_1 \prec P_2 \prec P_3 \prec P_4 \prec P_5,$
\end{center}

\noindent where the notation $P_i \prec P_{i+1}$ indicates that $P_{i+1}$ is infinitely near to $P_i$, that is, each $P_i$ belongs to the $i$-th infinitesimal neighborhood of $O$.

Indeed, by blowing up five times consecutively over $O$, we can check that the shared tangent directions at the infinitely near points are exactly the tangent directions to the strict transforms of $C$ at them. Moreover,  the infinitely near base points $P_1 \prec P_2 \prec P_3 \prec P_4 \prec P_5$ over $O$ are independent of $P$. 

We may ask ourselves if $\varphi_{Q\oplus P} = \varphi_Q \circ \varphi_P$, for all $P,Q \in C \setminus \{O\}$. If this relation is confirmed, it would yield a splitting as groups of the canonical complex of the pair $(\p^2,C)$ at $C$, coming from a set-theoretical section $\eta \colon C \hookrightarrow  \Dec(C)$. In general, the splitting property of a given short exact sequence can be relative to some subgroup and not global. This is the reason of the expression ``splitting at''. In our context, this subgroup can be continuous or discrete.

Let us investigate the possible set-theoretical section $\eta \colon C \hookrightarrow  \Dec(C)$ given by 
\begin{center}
    $\eta(P)=$
$\begin{cases}
\varphi_P, & \text{if $P \neq O$}\\
\Id_{\p^2}, & \text{otherwise}
\end{cases}$.
\end{center}

If $\eta$ is also a group homomorphism, we would have $\varphi_P^{-1}=\varphi_{\ominus P}$ for all $P \in C$, where $\ominus P$ denotes the inverse of $P$ under the group law $\oplus$ of $C$.

For all $P,Q \in C \setminus \{O\} $ with $P \neq \ominus Q$, let us compare the degrees of $\varphi_{Q\oplus P}$ and $\varphi_Q \circ \varphi_P$. We already know that $\deg(\varphi_{Q\oplus P})=4$. The following result will allow us to compute the second degree.

\begin{cor}[cf. \cite{alb} Corollary 4.2.12] Let $f$ be a plane Cremona map of homaloidal type $(d;m_1,\ldots,m_r)$, and let $g$ be a plane Cremona map of homaloidal type $(e;\ell_1,\ldots,\ell_s)$ satisfying that the first $k$ base points of $f$ coincide with those of $g$ and no further coincidence. Then the composite map $g \circ f^{-1}$ has degree
\begin{center}
    $de - \displaystyle \sum_{i=1}^k m_i\ell_i$.
\end{center}
\end{cor}

Since $\uBs(\varphi_{\ominus P}) = \{ \ominus P, O \prec P_1 \prec P_2 \prec P_3 \prec P_4 \prec P_5\}$, we have
\begin{center}
    $\uBs(\varphi_P) \cap \uBs(\varphi_{\ominus Q})=\{O \prec P_1 \prec P_2 \prec P_3 \prec P_4 \prec P_5\}$.
\end{center} 

Therefore, the previous result tells us that 
\begin{center}
    $\deg(\varphi_Q \circ \varphi_P)= 4\cdot 4 - 6 \cdot 1 \cdot 1 =10 \neq 4 = \deg(\varphi_{Q\oplus P})$,
\end{center}

\nin which implies that the maps $\varphi_{Q\oplus P}$ and $\varphi_Q \circ \varphi_P$ are distinct and do not make $\eta$ a group homomorphism.

Thus, our candidate $\eta$ does not yield a splitting of the canonical complex of $(\p^2,C)$ at $C$. We therefore conclude that we have $\varphi_Q \circ \varphi_P, \varphi_{Q\oplus P} \in \rho^{-1}(T_{Q \oplus P})$ with $\varphi_Q \circ \varphi_P \neq \varphi_{Q\oplus P}$. This allows us to produce many elements in $\Ine(C)$. For instance, given $P,Q,R \in C$ such that $P \oplus Q \oplus R = O$, then $\varphi_{P \oplus Q} \circ \varphi_R$ and $\varphi_P \circ \varphi_{Q \oplus R}$ belong to $\Ine(C)$. 

In \cite{bf}, Blanc \& Furter examined topologies and structures of the Cremona groups. 

\begin{defn}[cf. \cite{bf} Definition 2.1]\label{morfbir}
    Let $A$ and $X$ be irreducible algebraic varieties, and let $f$ be an $A$-birational self-map of the $A$-variety $A \times X$ which satisfy the following:

    \begin{enumerate}
        \item\label{morfbir1} $f$ induces an isomorphism $U \overset{\simeq} \longrightarrow V$, where $U$ and $V$ are open subsets of $A \times X$, whose projections on $A$ are surjective,

        \item\label{morfbir2} $f(a,x)=(a,\pr_2(f(a,x)))$, where $\pr_2$ denotes the second projection. Hence for each point $a \in A$, the birational map $x \dashmapsto \pr_2(f(a,x))$ corresponds to an element $f_a \in \Bir(X)$.
    \end{enumerate}

The map $a \mapsto f_a$ represents a map from $A$ to $\Bir(X)$ and it is called a \textit{morphism} from $A$ to $\Bir(X)$.
\end{defn}

These notions yield the following topology on $\Bir(X)$ called the \textit{Zariski topology}: a subset $F \subset \Bir(X)$ is closed in this topology if for any algebraic variety $A$ and any morphism $A \rightarrow \Bir(X)$, its preimage is closed.

In \cite{bf}, Blanc \& Furter studied the case where $X=\p^n$. Very recently and using more tools, Hassanzadeh \& Mostafazadehfard \cite{hm2} investigated similar aspects of $\Bir(X)$ when $X$ is an arbitrary projective variety over an infinite field $k$, of any characteristic and not necessarily algebraically closed.

We now present the following result, which confirms the expectation stated in \cite{bpv1}:

\begin{thm}\label{non split can comp}
    The canonical complex \ref{cc} of the pair $(\p^2,C)$ does not admit any splitting as groups at $C$ under the isomorphism $\Aut(C) \simeq C \rtimes \mathbb{Z}_d$.
    
\end{thm}

\begin{proof}
    For the sake of contradiction, suppose that we have a splitting given by a section $\eta \colon C \hookrightarrow \Dec(C)$. Such a map induces a $C$-birational self-map $f$ of the $C$-variety $C \times \p^2$ in the following way 
\begin{align*}
f \colon C \times \p^2 & \dash C \times \p^2 \\
(P,x)&\dashmapsto (P,\eta_P(x)),
\end{align*}
where $\eta_P \coloneqq \eta(P)$, for all $P \in C$.

Indeed, for all $P \in C$, the map $\eta_P$ is birational. Since $\Bs(\eta_P) \subset C$ for all $P \in C$ by Theorem \ref{pan thm 1.3}, $f$ determines an isomorphism of $U = V = C \times (\p^2 \setminus C)$ onto $C \times (\p^2 \setminus C)$ and therefore satisfies item \ref{morfbir1} of Definition \ref{morfbir}. 

This implies that $\eta$ is also a morphism from $C$ to $\Bir(\p^2)$ with respect to the Zariski topology, whose image is contained in $\Dec(C)$. By \cite[Lemma 2.19]{bf}, the image $\eta(C)$ of $\eta$ is a closed subgroup of $\Bir(\p^2)$, which has bounded degree. Thus $C \simeq \eta(C)$ as algebraic varieties, and therefore $\eta(C)$ is a projective algebraic group inside $\Bir(\p^2)$. However, this violates the fact that any algebraic subgroup of $\Bir(\p^2)$ is affine \cite[Remark 2.21]{bf}. Hence, there does not exist any section $\eta$, which implies the result.
\end{proof}

\section{Volume preserving versus standard Sarkisov factorization}\label{vp sark x sark fac}

So far this paper addresses a $2$-dimensional scenario. The following questions are generalizations of Theorem \ref{pan thm 1.3} and Theorem \ref{thm sark=vp cubic} in higher dimension, namely,
\begin{center}
    \begin{enumerate}
        \item\label{quest 1 chapter 5} \textit{Let $D \subset \p^n$ be a hypersurface of degree $n+1$ with canonical singularities and consider $\phi \in \Dec(D) \setminus \Aut(\p^n,D)$. Does it hold that $\Bs(\phi) \subset D$?}
        \item\label{quest 2 chapter 5} \textit{In dimension 3 and under the same assumptions, is the Sarkisov algorithm applied to $\phi$ automatically volume preserving? }
    \end{enumerate}
\end{center}

In this section, we present in detail the following example, which yields negative answers to the questions above.

\begin{ex}
    Let $(\p^3,D)$ be a Calabi-Yau pair, where $D$ is a general irreducible normal quartic surface having a single canonical singularity of type $A_1$. There exists $\phi \in \Dec(D)$ such that $\Bs(\phi) \not\subset D$ for which a Sarkisov decomposition obtained algorithmically is not volume preserving.
\end{ex}

Regarding Questions \ref{quest 1 chapter 5} and \ref{quest 2 chapter 5} in dimension $3$, we are dealing with a quartic surface $D \subset \p^3$. If $D$ is nonsingular, then it is a K3 surface and $\Bir(D)=\Aut(D)$. \cite{og, pq} produce examples of such quartic surfaces for which no nontrivial automorphism is derived from $\Bir(\p^3)$ by restriction. In these examples, we have $\Dec(D) = \{ \Id_{\p^3} \}$. Thus, neither of those questions is meaningful in this setting. We remark that in this case, $(\p^3,D)$ is a (t,c) Calabi-Yau pair. The situation changes if we allow strictly canonical singularities on $D$. \footnote{By strictly canonical, it is meant canonical but not terminal.} Recall that in this case, by Proposition \ref{vp can pair}, we have $\Bir^{\vp}(\p^3,D)=\Dec(D)$.

\paragraph{Construction of the example.} Based on the minimal resolution process, the simplest canonical surface singularity is of type $A_1$ \cite[Theorem 4.22]{km}. Consider $D \subset \p^3$ a general irreducible normal quartic surface having such a type of singularity at $P = (1:0:0:0)$. After suitable coordinate change, one can show that the equation of $D$ is of the form $x_0^2A+x_0B+C$, where $A,B,C \in \C[x_1,x_2,x_3]$ are general homogeneous polynomials of degrees 2, 3, 4, respectively. Moreover, $A$ is a quadratic form of rank 3.

In the proof of \cite[Claim 5.8]{acm}, Araujo, Corti \& Massarenti show that the birational involution 
\begin{center}
    $\phi \colon (x_0 : x_1 : x_2 : x_3) \mapsto (-Ax_0 - B : Ax_1 : Ax_2 : Ax_3)$
\end{center}

\noindent belongs to $\Dec(D)$. One can check that $\Bs(\phi)=V(A,B)$ and it consists of the union of six pairwise distinct lines through $P$ if we take $B$ general enough. This implies that $\Bs(\phi) \not\subset D$ as $D$ does not contain lines. This is in contrast with Theorem \ref{pan thm 1.3} in dimension 2, which asserts that the base locus is contained in the boundary divisor. This fact will allow us to construct a Sarkisov factorization that is not volume preserving, which shows that a generalization of the Theorem \ref{thm sark=vp cubic} does not hold in higher dimensions. Thus, the answer to both initial questions in this section is no.

In what follows, let us analyze the map $\phi$ carefully and obtain algorithmically two distinct Sarkisov decompositions: one volume preserving and the other non-volume preserving.

\paragraph{Analysis of the map $\phi$.} Its associated linear system $\Gamma \subset |\oo_{\p^3}(3)|$ is in particular contained in the linear system of space cubics passing through the reducible curve $V(A,B)$. Write $V(A,B) = L_1 \cup \ldots \cup L_6$, where each $L_i$ stands for a line through $P$.

Consider $\h$ a general member of $\Gamma$. Notice that $P$ is a singularity of $\Bs(\phi)$ as well as of $\h$ and $D$. In particular, we have that $P \in \Sing(\h)$ is a canonical singularity of type $A_1$ and $m_P(\h)=2$. This observation will be important later on on many occasions.

Indeed, $\h$ is of the form
\begin{center}
    $\lambda_0(-Ax_0-B) + \lambda_1Ax_1 + \lambda_2Ax_2 + \lambda_3Ax_3$,
\end{center}
for some $(\lambda_0 : \lambda_1 : \lambda_2 : \lambda_3) \in \p^3$ identified with $\Gamma$.

Dehomogenizing $\h$ with respect to $x_0$, we get the equation of $\h$ in $\{x_0 \neq 0\}$ becomes
\begin{center}
    $\lambda_0(-A-B) + \lambda_1Ax_1 + \lambda_2Ax_2 + \lambda_3Ax_3$.
\end{center}

Thus, $TC_P\h=\{-\lambda_0A=0\}$ whose projectivization is an irreducible conic, since $\rk(A)=3$. We can assure that $P \in \Sing(\h)$ is of type $A_1$ because a single blowup of the ambient space will be enough to resolve the singularity. One can check that the corresponding exceptional divisor intersected with the strict transform of $\h$ is an irreducible conic.

Another way to argue why $P \in \Sing(\h)$ is of type $A_1$ is by looking at $TC_P\h$ and comparing it with the tangent cones of normal forms of surface canonical singularities. We would be using the fact that canonical singularities are equivalent to rational double points in dimension 2.

Let us run the Sarkisov program for $\phi$. From now on we will follow the notation and algorithm described in \cite{cor1}. Although it has a slightly different notation, we refer the reader to \cite[Flowchart 13-1-9]{mat} for an explicit flowchart. 

\begin{nota}
    Henceforth, abusing notation, sometimes we will denote divisors on varieties and their strict transforms or pushforwards in others with the same symbol. Moreover, to avoid confusion in some instances, we will denote certain strict transforms or pushforwards with a right lower index indicating the ambient variety. We will do the same for general members of linear systems.
\end{nota}

We will exhibit in detail possible Sarkisov factorizations for $\phi$ proceeding in stages. Similarly to the surface case, there also exists a notion of \textit{Sarkisov degree} to detect the complexity of a birational map between threefolds with the structure of Mori fibered space. The point here is that it consists now of a triple of values and not a single one as in Definition \ref{defn sark deg dim 2}. We will briefly expose this notion, referring the reader to \cite{cor1,mat} for more explicit and precise definitions. 

\begin{defn}[cf. \cite{cor1} Definition 5.1]\label{defn sark deg dim 3}
    Let 
\[\begin{tikzcd}
	X & {X'} \\
	S & {S'}
	\arrow["f"', from=1-1, to=2-1]
	\arrow["\phi", dashed, from=1-1, to=1-2]
	\arrow["{f'}", from=1-2, to=2-2]
\end{tikzcd}\]

\nin be a birational map between threefold Mori fibered spaces. Consider the choice of $\mathcal{H}$ and $\mathcal{H'}$ made in \cite[Section 4]{cor1}. The \textit{Sarkisov degree} of $(\mathcal{H},~f \colon X \rightarrow S)$ is the triple $(\mu,c,e)$ where

\begin{enumerate}
    \item $\mu$ is the \textit{quasi-effective threshold} defined by $\mathcal{H} \equiv -\mu K_X$ over $S$ as in \cite[Section 4]{cor1};
    \item $c$ is the \textit{canonical threshold} of the pair $(X,\mathcal{H})$;
    \item $e$ is the number of \textit{crepant exceptional divisors} with respect to the pair $(X,c\mathcal{H})$.
\end{enumerate}
\end{defn}

If $X' \overset{f'} \longrightarrow S'$ is $\p^3 \rightarrow \Spec(\C)$, then $\h$ can be seen as a general member of the linear system associated to $\phi \colon X \dash \p^3$.

On the set of triples $(\mu,c,e)$ we introduce a \textit{partial ordering} as follows: 
\begin{center}
    $(\mu,c,e) > (\mu_1,c_1,e_1)$ if either
\end{center}
\begin{enumerate}
    \item $\mu > \mu_1$, or
    \item $\mu=\mu_1$ and $c < c_1$, or
    \item $\mu=\mu_1, c=c_1$ and $e>e_1$.
\end{enumerate}

\paragraph{Volume preserving factorization of $\phi$.} The starting point or Stage 0 in the Sarkisov program is to compute the Sarkisov degree $(\mu,c,e)$ of the corresponding birational map $\phi$. It will guide us along the factorization process. 

Let us compute $\mu$. By means of the projection formula, we need to solve for $\mu$ the equation $(\mathcal{H} + \mu K_{\p^3}) \cdot F = 0$, where $F$ is any curve contracted by the structure morphism $\p^3 \rightarrow \Spec(\C)$.  Therefore, $$\mu = - \dfrac{\h \cdot F}{K_{\p^3} \cdot F} =\dfrac{3\deg(F)}{4\deg(F)}=\dfrac{3}{4}.$$

By definition, $c= \max\{\lambda~ |~(\p^3,\lambda\h)~\text{is canonical}\}$. Performing a lot of computations and analyzing carefully the general member of the strict transform of $\h$, the map $\sigma \colon Z \rightarrow \p^3$ given by the consecutive blowups of

\begin{itemize}
    \item $P$,
    \item a curve $e$ in the exceptional divisor $E_P$ over $P$,
    \item a curve $e'$ in the exceptional divisor $E_e$ over $e$ and
    \item of the six lines $L_1,\ldots,L_6$
\end{itemize}

\noindent resolves the pair $(\p^3,\h)$. Set $E_{e'}$ to be the exceptional divisor over $e'$ and $E_i$ to be the exceptional divisor over $L_i$ for $i \in \{1,\ldots,6\}$. One has
\begin{center}
    $K_Z+\sigma^{-1}_*(\lambda\h)=\sigma^{*}(K_{\p^3}+\lambda\h)+(2-2\lambda)E_P+(3-3\lambda)E_e+(4-4\lambda)E_{e'}+(1-\lambda)E_1+\cdots+(1-\lambda)E_6,$
\end{center}

\noindent which easily yields $c=1$. 

Observing that $\sigma$ is also a \textit{ maximum crepant blow-up} for the pair $(\p^3,\h)$ \cite[Proposition-Definition 2.8]{cor1}, we deduce that the respective crepant exceptional divisors correspond to the exceptional divisors of $\sigma$. This gives us $e=9$. 

Thus, the Sarkisov degree of $\phi$ is $\left( \dfrac{3}{4},1,9 \right)$.

\subparagraph{Stage 1:} Following the Sarkisov algorithm in \cite{cor1}, since $c=1 < \dfrac{4}{3}=\dfrac{1}{\mu}$, the first link in the Sarkisov factorization is of type I or II. This link is initiated by an \textit{extremal blowup} \cite[Proposition-Definition 2.10]{cor1}, which always exists in this situation.

This choice of the extremal blowup is not determined by the algorithm. We are free to choose it. In our case, we have seven possibilities for such maps which are the blowup of $\p^3$ at $P$ or the blowup of $\p^3$ along one of the lines $L_i$ through $P$.


Only the first option gives us a volume preserving map, whereas the remaining ones do not. Indeed, let $\sigma_P \colon Z_1 \rightarrow \p^3$ be the blowup of $P$ and set $E_P \coloneqq \Exc(\sigma_P)$. By the adjunction formula, we have $K_{Z_1}=\sigma_P^*K_{\p^3}+2E_P$ and because $m_P(D)=2$, we have $D_{Z_1} \sim \sigma_P^*D-2E_P$. Hence,
\begin{center}
    $K_{Z_1}+D_{Z_1}=\sigma_P^*(K_{\p^3}+D)$,
\end{center}
\nin and since $(\sigma_P)_*D_{Z_1}=D$, Proposition-Definition \ref{prop crep bir morp} implies that $\sigma_P$ is volume preserving.

Without less of generality, let $\sigma_1 \colon Z_1' \rightarrow \p^3$ be the blowup of $L_1$ and set $E_1 \coloneqq \Exc(\sigma_1)$. By the adjunction formula, we have $K_{Z_1'}=\sigma_1^*K_{\p^3}+E_1$ and because $L_1 \not \subset D$, we have $D_{Z_1'} \sim \sigma_1^*D$. Hence,
\begin{center}
    $K_{Z_1'}+D_{Z_1'}=\sigma_1^*(K_{\p^3}+D)+E_1$,
\end{center}
showing that the crepant pullback contains a divisor with negative coefficient.

Let us proceed with $\sigma_P$. Consider $\h$ a general member of the linear system associated to $\phi$. Since $c=1$, the next thing to do is to run the $(K_{Z_1}+\h_{Z_1})$-MMP over $\Spec(\C)$. One can verify it results in the Mori fibered space structure $\pi_1 \colon Z_1 \rightarrow \p^2$. 

Thus, the first (volume preserving) Sarkisov link in a factorization of $\phi$ is of type I, and it is given by the blowup of $\p^3$ at $P$. We have that $E_P \simeq \p^2$, $\Tilde{A} \simeq \F_2$ and (by abuse of notation) $e = E_P \cap D_{Z_1} \simeq \p^1$. Moreover, the curve $e$ is contained in $\Bs(\phi \circ \sigma_P)$. All the lines $L_i$ are separated in $Z_1$ and are, in particular, rulings of $\F_2$. Observe that $\pi_1$ maps $e$ isomorphically onto the conic $\{A = 0\} \subset \p^2$. We have the following picture:


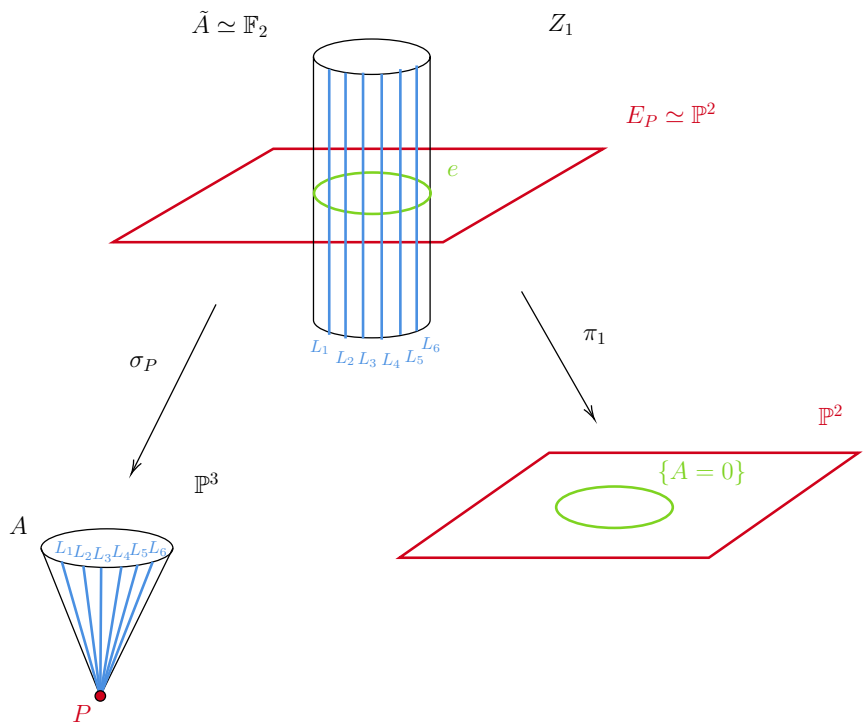
\begin{figure}[htb]
\centering
\resizebox{0.7\textwidth}{!}%
{

\tikzset{every picture/.style={line width=0.75pt}} 

\begin{tikzpicture}[x=0.75pt,y=0.75pt,yscale=-1,xscale=1]

\draw   (62,433.51) .. controls (62,425.28) and (84.61,418.6) .. (112.5,418.6) .. controls (140.39,418.6) and (163,425.28) .. (163,433.51) .. controls (163,441.75) and (140.39,448.42) .. (112.5,448.42) .. controls (84.61,448.42) and (62,441.75) .. (62,433.51) -- cycle ;
\draw [color={rgb, 255:red, 0; green, 0; blue, 0 }  ,draw opacity=1 ][line width=0.75]    (62,433.51) -- (107.39,547.6) -- (163,433.51) ;
\draw [color={rgb, 255:red, 74; green, 144; blue, 226 }  ,draw opacity=1 ][line width=1.5]    (77.86,443.77) -- (107.39,547.6) ;
\draw [color={rgb, 255:red, 74; green, 144; blue, 226 }  ,draw opacity=1 ][line width=1.5]    (94.44,447.28) -- (107.39,547.6) ;
\draw [color={rgb, 255:red, 74; green, 144; blue, 226 }  ,draw opacity=1 ][line width=1.5]    (123.57,447.74) -- (107.39,547.6) ;
\draw [color={rgb, 255:red, 74; green, 144; blue, 226 }  ,draw opacity=1 ][line width=1.5]    (147.53,444.09) -- (107.39,547.6) ;
\draw [color={rgb, 255:red, 74; green, 144; blue, 226 }  ,draw opacity=1 ][line width=1.5]    (108.03,448.65) -- (107.39,547.6) ;
\draw [color={rgb, 255:red, 74; green, 144; blue, 226 }  ,draw opacity=1 ][line width=1.5]    (135.87,446.37) -- (107.39,547.6) ;
\draw  [color={rgb, 255:red, 208; green, 2; blue, 27 }  ,draw opacity=1 ][line width=1.5]  (239.86,125.33) -- (492.63,125.33) -- (369.9,197.39) -- (117.13,197.39) -- cycle ;
\draw   (359.96,54.26) -- (359.96,257.55) .. controls (359.96,265.1) and (339.97,271.21) .. (315.3,271.21) .. controls (290.64,271.21) and (270.65,265.1) .. (270.65,257.55) -- (270.65,54.26) .. controls (270.65,46.72) and (290.64,40.6) .. (315.3,40.6) .. controls (339.97,40.6) and (359.96,46.72) .. (359.96,54.26) .. controls (359.96,61.81) and (339.97,67.93) .. (315.3,67.93) .. controls (290.64,67.93) and (270.65,61.81) .. (270.65,54.26) ;
\draw  [color={rgb, 255:red, 126; green, 211; blue, 33 }  ,draw opacity=1 ][line width=1.5]  (271.1,159.6) .. controls (271.1,150.78) and (291.09,143.62) .. (315.75,143.62) .. controls (340.41,143.62) and (360.41,150.78) .. (360.41,159.6) .. controls (360.41,168.43) and (340.41,175.59) .. (315.75,175.59) .. controls (291.09,175.59) and (271.1,168.43) .. (271.1,159.6) -- cycle ;
\draw [color={rgb, 255:red, 74; green, 144; blue, 226 }  ,draw opacity=1 ][line width=1.5]    (282.71,63.79) -- (282.71,268.9) ;
\draw [color={rgb, 255:red, 74; green, 144; blue, 226 }  ,draw opacity=1 ][line width=1.5]    (295.21,66.56) -- (295.21,271.68) ;
\draw [color={rgb, 255:red, 74; green, 144; blue, 226 }  ,draw opacity=1 ][line width=1.5]    (308.61,66.56) -- (308.61,271.68) ;
\draw [color={rgb, 255:red, 74; green, 144; blue, 226 }  ,draw opacity=1 ][line width=1.5]    (322.9,67.49) -- (322.9,272.6) ;
\draw [color={rgb, 255:red, 74; green, 144; blue, 226 }  ,draw opacity=1 ][line width=1.5]    (337.19,63.79) -- (337.19,268.9) ;
\draw [color={rgb, 255:red, 74; green, 144; blue, 226 }  ,draw opacity=1 ][line width=1.5]    (349.69,61.02) -- (349.69,266.13) ;
\draw  [color={rgb, 255:red, 208; green, 2; blue, 27 }  ,draw opacity=1 ][line width=1.5]  (451.17,359.93) -- (688.14,359.93) -- (573.46,440.93) -- (336.5,440.93) -- cycle ;
\draw  [color={rgb, 255:red, 126; green, 211; blue, 33 }  ,draw opacity=1 ][line width=1.5]  (456.51,401.93) .. controls (456.51,393.11) and (476.51,385.95) .. (501.17,385.95) .. controls (525.83,385.95) and (545.83,393.11) .. (545.83,401.93) .. controls (545.83,410.76) and (525.83,417.92) .. (501.17,417.92) .. controls (476.51,417.92) and (456.51,410.76) .. (456.51,401.93) -- cycle ;
\draw    (196,245.4) -- (131.89,374.61) ;
\draw [shift={(131,376.4)}, rotate = 296.39] [color={rgb, 255:red, 0; green, 0; blue, 0 }  ][line width=0.75]    (10.93,-3.29) .. controls (6.95,-1.4) and (3.31,-0.3) .. (0,0) .. controls (3.31,0.3) and (6.95,1.4) .. (10.93,3.29)   ;
\draw    (430,235.4) -- (485.02,332.66) ;
\draw [shift={(486,334.4)}, rotate = 240.51] [color={rgb, 255:red, 0; green, 0; blue, 0 }  ][line width=0.75]    (10.93,-3.29) .. controls (6.95,-1.4) and (3.31,-0.3) .. (0,0) .. controls (3.31,0.3) and (6.95,1.4) .. (10.93,3.29)   ;
\draw  [fill={rgb, 255:red, 208; green, 2; blue, 27 }  ,fill opacity=1 ] (103.39,547.6) .. controls (103.39,545.39) and (105.18,543.6) .. (107.39,543.6) .. controls (109.59,543.6) and (111.39,545.39) .. (111.39,547.6) .. controls (111.39,549.81) and (109.59,551.6) .. (107.39,551.6) .. controls (105.18,551.6) and (103.39,549.81) .. (103.39,547.6) -- cycle ;

\draw (84.4,553.28) node [anchor=north west][inner sep=0.75pt]  [font=\Large,color={rgb, 255:red, 208; green, 2; blue, 27 }  ,opacity=1 ]  {$P$};
\draw (36,407.78) node [anchor=north west][inner sep=0.75pt]  [font=\Large]  {$A$};
\draw (175.8,17.6) node [anchor=north west][inner sep=0.75pt]  [font=\Large]  {$\tilde{A} \simeq \mathbb{F}_{2}$};
\draw (508.13,88.87) node [anchor=north west][inner sep=0.75pt]  [font=\Large,color={rgb, 255:red, 208; green, 2; blue, 27 }  ,opacity=1 ]  {$E_{P} \simeq \mathbb{P}^{2}$};
\draw (371.6,136.6) node [anchor=north west][inner sep=0.75pt]  [font=\Large,color={rgb, 255:red, 126; green, 211; blue, 33 }  ,opacity=1 ]  {$e$};
\draw (532.6,364.53) node [anchor=north west][inner sep=0.75pt]  [font=\Large,color={rgb, 255:red, 126; green, 211; blue, 33 }  ,opacity=1 ]  {$\{A=0\}$};
\draw (655.87,322) node [anchor=north west][inner sep=0.75pt]  [font=\Large,color={rgb, 255:red, 208; green, 2; blue, 27 }  ,opacity=1 ]  {$\mathbb{P}^{2}$};
\draw (128.2,284.2) node [anchor=north west][inner sep=0.75pt]  [font=\Large]  {$\sigma _{P}$};
\draw (475.67,263.4) node [anchor=north west][inner sep=0.75pt]  [font=\Large]  {$\pi _{1}$};
\draw (70.83,426.4) node [anchor=north west][inner sep=0.75pt]  [font=\small,color={rgb, 255:red, 74; green, 144; blue, 226 }  ,opacity=1 ]  {$L_{1}$};
\draw (85,429.9) node [anchor=north west][inner sep=0.75pt]  [font=\small,color={rgb, 255:red, 74; green, 144; blue, 226 }  ,opacity=1 ]  {$L_{2}$};
\draw (100,431.4) node [anchor=north west][inner sep=0.75pt]  [font=\small,color={rgb, 255:red, 74; green, 144; blue, 226 }  ,opacity=1 ]  {$L_{3}$};
\draw (114.5,429.9) node [anchor=north west][inner sep=0.75pt]  [font=\small,color={rgb, 255:red, 74; green, 144; blue, 226 }  ,opacity=1 ]  {$L_{4}$};
\draw (128,427.9) node [anchor=north west][inner sep=0.75pt]  [font=\small,color={rgb, 255:red, 74; green, 144; blue, 226 }  ,opacity=1 ]  {$L_{5}$};
\draw (142.33,426.73) node [anchor=north west][inner sep=0.75pt]  [font=\small,color={rgb, 255:red, 74; green, 144; blue, 226 }  ,opacity=1 ]  {$L_{6}$};
\draw (266.5,271.4) node [anchor=north west][inner sep=0.75pt]  [font=\small,color={rgb, 255:red, 74; green, 144; blue, 226 }  ,opacity=1 ]  {$L_{1}$};
\draw (286,279.9) node [anchor=north west][inner sep=0.75pt]  [font=\small,color={rgb, 255:red, 74; green, 144; blue, 226 }  ,opacity=1 ]  {$L_{2}$};
\draw (303,280.73) node [anchor=north west][inner sep=0.75pt]  [font=\small,color={rgb, 255:red, 74; green, 144; blue, 226 }  ,opacity=1 ]  {$L_{3}$};
\draw (321.17,281.57) node [anchor=north west][inner sep=0.75pt]  [font=\small,color={rgb, 255:red, 74; green, 144; blue, 226 }  ,opacity=1 ]  {$L_{4}$};
\draw (339.08,278.97) node [anchor=north west][inner sep=0.75pt]  [font=\small,color={rgb, 255:red, 74; green, 144; blue, 226 }  ,opacity=1 ,xslant=-0.01]  {$L_{5}$};
\draw (351.69,269.53) node [anchor=north west][inner sep=0.75pt]  [font=\small,color={rgb, 255:red, 74; green, 144; blue, 226 }  ,opacity=1 ]  {$L_{6}$};
\draw (448.67,19.73) node [anchor=north west][inner sep=0.75pt]  [font=\Large]  {$Z_{1}$};
\draw (178,374.07) node [anchor=north west][inner sep=0.75pt]  [font=\Large]  {$\mathbb{P}^{3}$};

\end{tikzpicture}
}  
\caption{First link in a Sarkisov factorization of $\phi$.}
\label{fig:first link}
\end{figure}

\subparagraph{Stage 2:} We must compute the Sarkisov degree $(\mu_1,c_1,e_1)$ of the induced birational map $\phi \circ \sigma_P$.  One can check that it equals $\left( \dfrac{1}{2},1,8  \right)$ and it is smaller than $(\mu,c,e)$ according to the partial ordering explained after the Definition \ref{defn sark deg dim 3}. So the birational map $\phi \circ \sigma_P$ is ``simpler'' than $\phi$.

Since $c_1 = 1 < 2 = \dfrac{1}{\mu_1}$, the second link in the Sarkisov factorization is of type I or II. At this point, we also have seven extremal blowups to choose from. They are the blowup of $Z_1$ along $e$ or the blowup of $Z_1$ along one of the lines $L_i$. Repeating exactly the same arguments as in Stage 1, we can check that the first one yields a volume preserving map whereas the remaining ones do not.

Let us continue with $\sigma_e \colon Y_1 \rightarrow Z_1$ the blowup of $Z_1$ along $e$ and set $E_e \coloneqq \Exc(\sigma_e)$. 

By \cite[Theorem 8.24, (b)]{har}, we have that $E_e \simeq \p(\n_{e/Z_1}^{\vee})$. One can compute $\n_{e/Z_1}^{\vee} \simeq \oo_{\p^1}(2) \oplus \oo_{\p^1}(-4)$ and therefore $E_e \simeq \F_6$. Indeed, this follows from the existence of the following short exact sequence of normal bundles and identifications $\n_{e/E_P}\simeq \oo_{\p^1}(4)$ and $(\n_{E_P/Z_1})|_{e}\simeq \oo_{\p^1}(-2)$,

\begin{equation*}
        0 \longrightarrow \n_{e/E_P} \longrightarrow \n_{e/Z_1} \longrightarrow (\n_{E_P/Z_1})|_{e} \longrightarrow 0.
\end{equation*}

The curve $e' = E_e \cap D_{Y_1} \simeq \p^1$ is contained in $\Bs(\phi \circ \sigma_P \circ \sigma_e)$.

Since $c_1=1$, we need to run the $(K_{Y_1}+\h_{Y_1})$-MMP over $\p^2$. One can verify that this log MMP results in the divisorial contraction 
\begin{center}
    $\alpha \colon Y_1 \rightarrow Z_2 \simeq \p(\oo_{\p^2} \oplus \oo_{\p^2}(3))$,
\end{center}
where $\pi_2 \colon Z_2 \rightarrow \p^2$ is the corresponding structure morphism making $Z_2$ a Mori fibered space. This birational morphism contracts exactly the rulings of $\Tilde{A} \simeq \F_2$, that is, $\Tilde{A} = \Exc(\alpha)$. The isomorphism $Z_2 \simeq \p(\oo_{\p^2} \oplus \oo_{\p^2}(3))$ may be justified by analyzing the section of $\pi_2$ given by $\p^2 \rightarrow E_P \simeq \p^2$; and by $K_{Z_2} = \alpha_*K_{Y_1}$ written in terms of a basis for $\Pic(Z_2)$ and comparing it with the formula for the canonical class of a projective bundle. 

Thus, the second (volume preserving) Sarkisov link in a factorization of $\phi$ is of type II, and it is given by the composition $\alpha \circ \sigma_e^{-1}$. Observe that $\alpha$ maps $E_e$ isomorphically, via pushforward, onto the cylinder $E_e = \{ A = 0 \} \subset Z_2$. In particular, all the lines $L_i$ are contracted by $\alpha$. Moreover, we have that $\Bs(\phi \circ \sigma_P \circ \sigma_e \circ \alpha^{-1})$ consists of the curve $f=\alpha(e')$ which is mapped by $\pi_2$ isomorphically onto the conic $\{A=0\} \subset \p^2$. We have the following picture in which we did not put the strict transforms of $D$:

\newpage

\begin{figure}[htb]
\centering
\resizebox{0.9\textwidth}{!}%
{

\tikzset{every picture/.style={line width=0.75pt}} 

\begin{tikzpicture}[x=0.75pt,y=0.75pt,yscale=-1,xscale=1]

\draw  [color={rgb, 255:red, 208; green, 2; blue, 27 }  ,draw opacity=1 ][line width=1.5]  (132.56,473.25) -- (356.15,473.25) -- (247.59,536.94) -- (24,536.94) -- cycle ;
\draw   (234.49,414.22) -- (234.49,593.85) .. controls (234.49,600.53) and (216.81,605.94) .. (194.99,605.94) .. controls (173.17,605.94) and (155.49,600.53) .. (155.49,593.85) -- (155.49,414.22) .. controls (155.49,407.55) and (173.17,402.13) .. (194.99,402.13) .. controls (216.81,402.13) and (234.49,407.55) .. (234.49,414.22) .. controls (234.49,420.9) and (216.81,426.31) .. (194.99,426.31) .. controls (173.17,426.31) and (155.49,420.9) .. (155.49,414.22) ;
\draw  [color={rgb, 255:red, 126; green, 211; blue, 33 }  ,draw opacity=1 ][line width=1.5]  (155.88,507.3) .. controls (155.88,499.5) and (173.57,493.18) .. (195.38,493.18) .. controls (217.2,493.18) and (234.89,499.5) .. (234.89,507.3) .. controls (234.89,515.1) and (217.2,521.43) .. (195.38,521.43) .. controls (173.57,521.43) and (155.88,515.1) .. (155.88,507.3) -- cycle ;
\draw [color={rgb, 255:red, 74; green, 144; blue, 226 }  ,draw opacity=1 ][line width=1.5]    (166.15,422.63) -- (166.15,603.9) ;
\draw [color={rgb, 255:red, 74; green, 144; blue, 226 }  ,draw opacity=1 ][line width=1.5]    (177.21,425.08) -- (177.21,606.35) ;
\draw [color={rgb, 255:red, 74; green, 144; blue, 226 }  ,draw opacity=1 ][line width=1.5]    (189.06,425.08) -- (189.06,606.35) ;
\draw [color={rgb, 255:red, 74; green, 144; blue, 226 }  ,draw opacity=1 ][line width=1.5]    (201.7,425.89) -- (201.7,607.17) ;
\draw [color={rgb, 255:red, 74; green, 144; blue, 226 }  ,draw opacity=1 ][line width=1.5]    (214.35,422.63) -- (214.35,603.9) ;
\draw [color={rgb, 255:red, 74; green, 144; blue, 226 }  ,draw opacity=1 ][line width=1.5]    (225.41,420.18) -- (225.41,601.45) ;
\draw  [color={rgb, 255:red, 208; green, 2; blue, 27 }  ,draw opacity=1 ][line width=1.5]  (424.56,230.66) -- (648.15,230.66) -- (539.59,304.54) -- (316,304.54) -- cycle ;
\draw   (397.38,44.47) -- (617.64,112.53) .. controls (625.14,114.84) and (623.89,136.14) .. (614.84,160.1) .. controls (605.79,184.05) and (592.38,201.59) .. (584.88,199.28) -- (364.62,131.22) .. controls (357.12,128.9) and (358.37,107.6) .. (367.42,83.65) .. controls (376.46,59.69) and (389.88,42.15) .. (397.38,44.47) .. controls (404.88,46.79) and (403.63,68.09) .. (394.58,92.04) .. controls (385.53,116) and (372.12,133.54) .. (364.62,131.22) ;
\draw  [color={rgb, 255:red, 144; green, 19; blue, 254 }  ,draw opacity=1 ][line width=1.5]  (455.67,158.61) .. controls (446.17,153.41) and (456.95,133.34) .. (479.74,113.79) .. controls (502.53,94.23) and (528.71,82.59) .. (538.2,87.79) .. controls (547.7,92.99) and (536.93,113.06) .. (514.13,132.62) .. controls (491.34,152.17) and (465.17,163.81) .. (455.67,158.61) -- cycle ;
\draw [color={rgb, 255:red, 74; green, 144; blue, 226 }  ,draw opacity=1 ][line width=1.5]    (378.65,121.44) -- (600.9,191) ;
\draw [color={rgb, 255:red, 74; green, 144; blue, 226 }  ,draw opacity=1 ][line width=1.5]    (385.15,111.44) -- (608.46,179.77) ;
\draw [color={rgb, 255:red, 74; green, 144; blue, 226 }  ,draw opacity=1 ][line width=1.5]    (392.15,98.76) -- (612.21,166.76) ;
\draw [color={rgb, 255:red, 74; green, 144; blue, 226 }  ,draw opacity=1 ][line width=1.5]    (397.65,83.94) -- (619.61,153.18) ;
\draw [color={rgb, 255:red, 74; green, 144; blue, 226 }  ,draw opacity=1 ][line width=1.5]    (400.82,70.08) -- (620.89,138.08) ;
\draw [color={rgb, 255:red, 74; green, 144; blue, 226 }  ,draw opacity=1 ][line width=1.5]    (402.44,57.02) -- (622.5,125.01) ;
\draw  [color={rgb, 255:red, 126; green, 211; blue, 33 }  ,draw opacity=1 ][line width=1.5]  (542.5,31.11) -- (542.5,266.66) .. controls (542.5,274.16) and (522.24,280.23) .. (497.25,280.23) .. controls (472.26,280.23) and (452,274.16) .. (452,266.66) -- (452,31.11) .. controls (452,23.62) and (472.26,17.54) .. (497.25,17.54) .. controls (522.24,17.54) and (542.5,23.62) .. (542.5,31.11) .. controls (542.5,38.61) and (522.24,44.69) .. (497.25,44.69) .. controls (472.26,44.69) and (452,38.61) .. (452,31.11) ;
\draw  [draw opacity=0][line width=1.5]  (542.1,264.38) .. controls (542.34,264.95) and (542.48,265.53) .. (542.49,266.13) .. controls (542.61,273.7) and (522.61,280.16) .. (497.83,280.56) .. controls (473.05,280.95) and (452.86,275.13) .. (452.74,267.55) .. controls (452.73,266.85) and (452.89,266.16) .. (453.21,265.49) -- (497.61,266.84) -- cycle ; \draw  [color={rgb, 255:red, 65; green, 117; blue, 5 }  ,draw opacity=1 ][line width=1.5]  (542.1,264.38) .. controls (542.34,264.95) and (542.48,265.53) .. (542.49,266.13) .. controls (542.61,273.7) and (522.61,280.16) .. (497.83,280.56) .. controls (473.05,280.95) and (452.86,275.13) .. (452.74,267.55) .. controls (452.73,266.85) and (452.89,266.16) .. (453.21,265.49) ;  
\draw  [draw opacity=0][dash pattern={on 5.63pt off 4.5pt}][line width=1.5]  (453.61,265.95) .. controls (453.37,265.38) and (453.25,264.79) .. (453.25,264.2) .. controls (453.29,256.62) and (473.41,250.57) .. (498.19,250.68) .. controls (522.98,250.79) and (543.04,257.03) .. (543.01,264.61) .. controls (543.01,265.31) and (542.83,265.99) .. (542.5,266.66) -- (498.13,264.4) -- cycle ; \draw  [color={rgb, 255:red, 65; green, 117; blue, 5 }  ,draw opacity=1 ][dash pattern={on 5.63pt off 4.5pt}][line width=1.5]  (453.61,265.95) .. controls (453.37,265.38) and (453.25,264.79) .. (453.25,264.2) .. controls (453.29,256.62) and (473.41,250.57) .. (498.19,250.68) .. controls (522.98,250.79) and (543.04,257.03) .. (543.01,264.61) .. controls (543.01,265.31) and (542.83,265.99) .. (542.5,266.66) ;  
\draw  [color={rgb, 255:red, 208; green, 2; blue, 27 }  ,draw opacity=1 ][line width=1.5]  (728.9,564.6) -- (931.77,564.6) -- (833.28,632.82) -- (630.41,632.82) -- cycle ;
\draw  [color={rgb, 255:red, 126; green, 211; blue, 33 }  ,draw opacity=1 ][line width=1.5]  (833.74,380.11) -- (833.74,598.06) .. controls (833.74,604.86) and (815.35,610.38) .. (792.68,610.38) .. controls (770,610.38) and (751.62,604.86) .. (751.62,598.06) -- (751.62,380.11) .. controls (751.62,373.31) and (770,367.8) .. (792.68,367.8) .. controls (815.35,367.8) and (833.74,373.31) .. (833.74,380.11) .. controls (833.74,386.92) and (815.35,392.43) .. (792.68,392.43) .. controls (770,392.43) and (751.62,386.92) .. (751.62,380.11) ;
\draw  [draw opacity=0][line width=1.5]  (833.38,595.76) .. controls (833.6,596.28) and (833.71,596.81) .. (833.72,597.34) .. controls (833.83,604.34) and (815.69,610.31) .. (793.21,610.67) .. controls (770.72,611.04) and (752.4,605.66) .. (752.3,598.67) .. controls (752.29,598.03) and (752.43,597.41) .. (752.7,596.79) -- (793.01,598.01) -- cycle ; \draw  [color={rgb, 255:red, 65; green, 117; blue, 5 }  ,draw opacity=1 ][line width=1.5]  (833.38,595.76) .. controls (833.6,596.28) and (833.71,596.81) .. (833.72,597.34) .. controls (833.83,604.34) and (815.69,610.31) .. (793.21,610.67) .. controls (770.72,611.04) and (752.4,605.66) .. (752.3,598.67) .. controls (752.29,598.03) and (752.43,597.41) .. (752.7,596.79) ;  
\draw  [draw opacity=0][dash pattern={on 5.63pt off 4.5pt}][line width=1.5]  (753.07,597.16) .. controls (752.86,596.64) and (752.76,596.11) .. (752.76,595.57) .. controls (752.79,588.58) and (771.05,582.99) .. (793.54,583.09) .. controls (816.02,583.19) and (834.23,588.95) .. (834.2,595.94) .. controls (834.19,596.58) and (834.04,597.2) .. (833.75,597.8) -- (793.48,595.76) -- cycle ; \draw  [color={rgb, 255:red, 65; green, 117; blue, 5 }  ,draw opacity=1 ][dash pattern={on 5.63pt off 4.5pt}][line width=1.5]  (753.07,597.16) .. controls (752.86,596.64) and (752.76,596.11) .. (752.76,595.57) .. controls (752.79,588.58) and (771.05,582.99) .. (793.54,583.09) .. controls (816.02,583.19) and (834.23,588.95) .. (834.2,595.94) .. controls (834.19,596.58) and (834.04,597.2) .. (833.75,597.8) ;  
\draw  [draw opacity=0][line width=1.5]  (834.25,466.4) .. controls (834.31,466.66) and (834.34,466.92) .. (834.34,467.18) .. controls (834.33,473.44) and (815.94,478.32) .. (793.27,478.09) .. controls (770.87,477.86) and (752.66,472.71) .. (752.26,466.55) -- (793.3,466.75) -- cycle ; \draw  [color={rgb, 255:red, 144; green, 19; blue, 254 }  ,draw opacity=1 ][line width=1.5]  (834.25,466.4) .. controls (834.31,466.66) and (834.34,466.92) .. (834.34,467.18) .. controls (834.33,473.44) and (815.94,478.32) .. (793.27,478.09) .. controls (770.87,477.86) and (752.66,472.71) .. (752.26,466.55) ;  
\draw  [draw opacity=0][dash pattern={on 5.63pt off 4.5pt}][line width=1.5]  (752.47,464.54) .. controls (754.2,457.01) and (771.94,451.51) .. (793.5,452.08) .. controls (813.88,452.61) and (830.72,458.39) .. (833.77,465.46) -- (793.3,466.75) -- cycle ; \draw  [color={rgb, 255:red, 144; green, 19; blue, 254 }  ,draw opacity=1 ][dash pattern={on 5.63pt off 4.5pt}][line width=1.5]  (752.47,464.54) .. controls (754.2,457.01) and (771.94,451.51) .. (793.5,452.08) .. controls (813.88,452.61) and (830.72,458.39) .. (833.77,465.46) ;  
\draw [color={rgb, 255:red, 245; green, 166; blue, 35 }  ,draw opacity=1 ][line width=1.5]    (751.31,482.42) .. controls (751.6,471.91) and (763,456.74) .. (775.8,455.14) .. controls (788.6,453.54) and (774.69,489.38) .. (784.67,486) .. controls (794.65,482.61) and (793.75,459.54) .. (801.47,454) .. controls (809.18,448.46) and (807.97,483.1) .. (813.47,482.8) .. controls (818.96,482.49) and (830.27,418.8) .. (833.77,465.46) .. controls (837.27,512.12) and (819.87,512.4) .. (792.67,508.4) .. controls (765.47,504.4) and (752.25,489.71) .. (751.31,482.42) -- cycle ;
\draw    (284.17,245.47) -- (218.03,355.69) ;
\draw [shift={(217,357.4)}, rotate = 300.97] [color={rgb, 255:red, 0; green, 0; blue, 0 }  ][line width=0.75]    (10.93,-3.29) .. controls (6.95,-1.4) and (3.31,-0.3) .. (0,0) .. controls (3.31,0.3) and (6.95,1.4) .. (10.93,3.29)   ;
\draw    (724,227.4) -- (779.02,324.66) ;
\draw [shift={(780,326.4)}, rotate = 240.51] [color={rgb, 255:red, 0; green, 0; blue, 0 }  ][line width=0.75]    (10.93,-3.29) .. controls (6.95,-1.4) and (3.31,-0.3) .. (0,0) .. controls (3.31,0.3) and (6.95,1.4) .. (10.93,3.29)   ;
\draw  [dash pattern={on 4.5pt off 4.5pt}]  (427,512.44) -- (606.17,512.58) ;
\draw [shift={(608.17,512.58)}, rotate = 180.04] [color={rgb, 255:red, 0; green, 0; blue, 0 }  ][line width=0.75]    (10.93,-3.29) .. controls (6.95,-1.4) and (3.31,-0.3) .. (0,0) .. controls (3.31,0.3) and (6.95,1.4) .. (10.93,3.29)   ;
\draw [color={rgb, 255:red, 245; green, 166; blue, 35 }  ,draw opacity=1 ][line width=1.5]    (453,153.41) .. controls (453.32,142.03) and (459.54,139.49) .. (463,140.66) .. controls (466.46,141.83) and (476.4,154.9) .. (482.07,154.47) .. controls (487.74,154.04) and (479.5,135.16) .. (489.25,131.91) .. controls (499,128.66) and (497.57,152.72) .. (510.07,150.47) .. controls (522.57,148.22) and (504.15,119.97) .. (514,115.16) .. controls (523.85,110.35) and (523.22,131.02) .. (530.73,130.47) .. controls (538.25,129.92) and (529.98,97.09) .. (534.29,93.08) .. controls (538.6,89.07) and (541.93,98.66) .. (542.98,110.43) .. controls (544.02,122.19) and (540.64,139.09) .. (526.64,151.43) .. controls (512.64,163.76) and (507.98,168.76) .. (480.98,171.43) .. controls (453.98,174.09) and (454.04,161.31) .. (453,153.41) -- cycle ;

\draw (74.04,372.27) node [anchor=north west][inner sep=0.75pt]  [font=\Large]  {$\tilde{A} \simeq \mathbb{F}_{2}$};
\draw (353.61,445.68) node [anchor=north west][inner sep=0.75pt]  [font=\Large,color={rgb, 255:red, 208; green, 2; blue, 27 }  ,opacity=1 ]  {$E_{P} \simeq \mathbb{P}^{2}$};
\draw (243.98,485.97) node [anchor=north west][inner sep=0.75pt]  [font=\Large,color={rgb, 255:red, 126; green, 211; blue, 33 }  ,opacity=1 ]  {$e$};
\draw (149.62,602.91) node [anchor=north west][inner sep=0.75pt]  [font=\small,color={rgb, 255:red, 74; green, 144; blue, 226 }  ,opacity=1 ]  {$L_{1}$};
\draw (165.86,612.94) node [anchor=north west][inner sep=0.75pt]  [font=\small,color={rgb, 255:red, 74; green, 144; blue, 226 }  ,opacity=1 ]  {$L_{2}$};
\draw (183.12,615.17) node [anchor=north west][inner sep=0.75pt]  [font=\small,color={rgb, 255:red, 74; green, 144; blue, 226 }  ,opacity=1 ]  {$L_{3}$};
\draw (200.5,615.77) node [anchor=north west][inner sep=0.75pt]  [font=\small,color={rgb, 255:red, 74; green, 144; blue, 226 }  ,opacity=1 ]  {$L_{4}$};
\draw (215.31,613.28) node [anchor=north west][inner sep=0.75pt]  [font=\small,color={rgb, 255:red, 74; green, 144; blue, 226 }  ,opacity=1 ]  {$L_{5}$};
\draw (227.41,604.85) node [anchor=north west][inner sep=0.75pt]  [font=\small,color={rgb, 255:red, 74; green, 144; blue, 226 }  ,opacity=1 ]  {$L_{6}$};
\draw (290.2,33.04) node [anchor=north west][inner sep=0.75pt]  [font=\Large]  {$\tilde{A} \simeq \mathbb{F}_{2}$};
\draw (630.82,195.06) node [anchor=north west][inner sep=0.75pt]  [font=\Large,color={rgb, 255:red, 208; green, 2; blue, 27 }  ,opacity=1 ]  {$E_{P} \simeq \mathbb{P}^{2}$};
\draw (302.48,163.57) node [anchor=north west][inner sep=0.75pt]  [font=\Large,color={rgb, 255:red, 245; green, 166; blue, 35 }  ,opacity=1 ]  {$E_{e} \ \cap D_{Y_{1}} =e'$};
\draw (600.65,192.64) node [anchor=north west][inner sep=0.75pt]  [font=\small,color={rgb, 255:red, 74; green, 144; blue, 226 }  ,opacity=1 ]  {$L_{1}$};
\draw (609.39,177.71) node [anchor=north west][inner sep=0.75pt]  [font=\small,color={rgb, 255:red, 74; green, 144; blue, 226 }  ,opacity=1 ]  {$L_{2}$};
\draw (616.84,163.5) node [anchor=north west][inner sep=0.75pt]  [font=\small,color={rgb, 255:red, 74; green, 144; blue, 226 }  ,opacity=1 ]  {$L_{3}$};
\draw (622.78,148.71) node [anchor=north west][inner sep=0.75pt]  [font=\small,color={rgb, 255:red, 74; green, 144; blue, 226 }  ,opacity=1 ]  {$L_{4}$};
\draw (626.81,134.62) node [anchor=north west][inner sep=0.75pt]  [font=\small,color={rgb, 255:red, 74; green, 144; blue, 226 }  ,opacity=1 ]  {$L_{5}$};
\draw (626.76,119.23) node [anchor=north west][inner sep=0.75pt]  [font=\small,color={rgb, 255:red, 74; green, 144; blue, 226 }  ,opacity=1 ]  {$L_{6}$};
\draw (560.15,17.06) node [anchor=north west][inner sep=0.75pt]  [font=\Large,color={rgb, 255:red, 126; green, 211; blue, 33 }  ,opacity=1 ]  {$E_{e} \simeq \mathbb{F}_{6}$};
\draw (919.47,575.31) node [anchor=north west][inner sep=0.75pt]  [font=\Large,color={rgb, 255:red, 208; green, 2; blue, 27 }  ,opacity=1 ]  {$E_{P} \simeq \mathbb{P}^{2}$};
\draw (642.24,478.06) node [anchor=north west][inner sep=0.75pt]  [font=\Large,color={rgb, 255:red, 245; green, 166; blue, 35 }  ,opacity=1 ]  {$f\ =\ \alpha ( e')$};
\draw (847.02,366.7) node [anchor=north west][inner sep=0.75pt]  [font=\Large,color={rgb, 255:red, 126; green, 211; blue, 33 }  ,opacity=1 ]  {$E_{e} \simeq \mathbb{F}_{6}$};
\draw (549.15,56.91) node [anchor=north west][inner sep=0.75pt]  [font=\Large,color={rgb, 255:red, 144; green, 19; blue, 254 }  ,opacity=1 ]  {$E_{e} \ \cap \ \tilde{A}$};
\draw (848.6,517.2) node [anchor=north west][inner sep=0.75pt]  [font=\Large,color={rgb, 255:red, 126; green, 211; blue, 33 }  ,opacity=1 ]  {$\{A=0\}$};
\draw (219.2,270.2) node [anchor=north west][inner sep=0.75pt]  [font=\Large]  {$\sigma _{e}$};
\draw (768.67,260.4) node [anchor=north west][inner sep=0.75pt]  [font=\Large]  {$\alpha $};
\draw (318,384.73) node [anchor=north west][inner sep=0.75pt]  [font=\Large]  {$Z_{1}$};
\draw (681,16.73) node [anchor=north west][inner sep=0.75pt]  [font=\Large]  {$Y_{1}$};
\draw (965,381.73) node [anchor=north west][inner sep=0.75pt]  [font=\Large]  {$Z_{2}$};

\end{tikzpicture}

}  
\caption{Second link in a Sarkisov factorization of $\phi$.}
\label{fig:second link vp}
\end{figure}
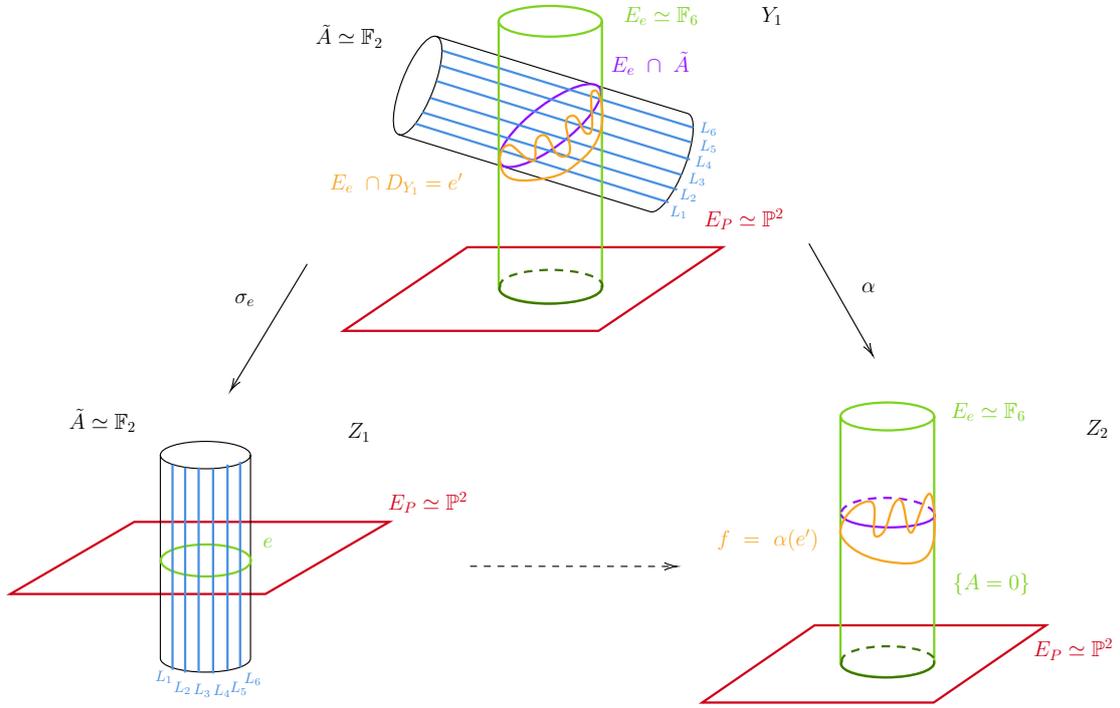

\subparagraph{Stage 3:} Once again, we need to compute the Sarkisov degree $(\mu_2,c_2,e_2)$ of the induced birational map $\phi \circ \sigma_P \circ \sigma_e \circ \alpha^{-1}$.  One can check that equals $\left( \dfrac{1}{2},1,1  \right)$ and it is smaller than $(\mu_1,c_1,e_1)$. Thus we have simplified the birational map $\phi \circ \sigma_P$. 

Since $c_2 = 1 < 2 = \dfrac{1}{\mu_2}$, the third link in the Sarkisov factorization is of type I or II. The difference in this stage is that we have only one possible extremal blowup given by the blowup of $Z_2$ along $f$. Consider $\sigma_f \colon Y_2 \rightarrow Z_2$ such map and denote $E_f \coloneqq \Exc(\sigma_f)$.

As in the previous stages, the map $\sigma_f$ is volume preserving because $f \subset D_{Z_2}$. One can check that the map $\phi \circ \sigma_P \circ \sigma_e \circ \alpha^{-1} \circ \sigma_f$ is everywhere defined. 

Since $c_1=1$, we need to run the $(K_{Y_2}+\h_{Y_2})$-MMP over $\p^2$. By \cite[Theorem 8.24, (b)]{har}, we have that $E_f \simeq \p(\n_{f/Z_2}^{\vee})$. One can compute $\n_{f/Z_2}^{\vee} \simeq \oo_{\p^1}(-6) \oplus \oo_{\p^1}(-4)$ and therefore $E_f \simeq \F_2$. Indeed, this follows from the existence of the following short exact sequence of normal bundles and identifications $\n_{f/E_e}\simeq \oo_{\p^1}(6)$ and $(\n_{E_e/Z_2})|_{f}\simeq \oo_{\p^1}(4)$,

\begin{equation*}
        0 \longrightarrow \n_{f/E_e} \longrightarrow \n_{f/Z_2} \longrightarrow (\n_{E_e/Z_2})|_{f} \longrightarrow 0.
\end{equation*}

Repeating the same arguments as in Stage 2, we obtain a divisorial contraction $\beta \colon Y_2 \rightarrow Z_3 \simeq \p(\oo_{\p^2} \oplus \oo_{\p^2}(1))$. The third (volume preserving) Sarkisov link in a factorization of $\phi$ is of type II, and it is given by the composition $\beta \circ \sigma_f^{-1}$ with analogous geometric properties to the previous one. We have the following picture:

\newpage

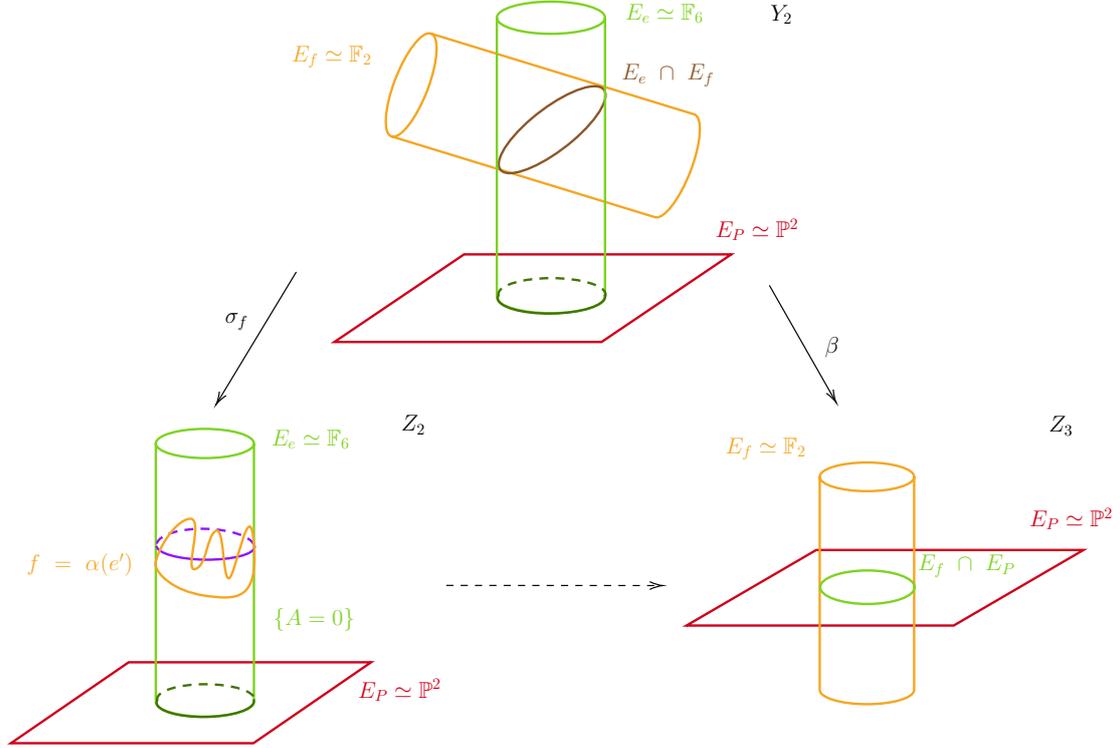
\begin{figure}[htb]
\centering
\resizebox{0.9\textwidth}{!}%
{

\tikzset{every picture/.style={line width=0.75pt}} 

\begin{tikzpicture}[x=0.75pt,y=0.75pt,yscale=-1,xscale=1]

\draw  [color={rgb, 255:red, 208; green, 2; blue, 27 }  ,draw opacity=1 ][line width=1.5]  (424.56,230.66) -- (648.15,230.66) -- (539.59,304.54) -- (316,304.54) -- cycle ;
\draw  [color={rgb, 255:red, 245; green, 166; blue, 35 }  ,draw opacity=1 ][line width=1.5]  (396.71,44.47) -- (616.97,112.53) .. controls (624.47,114.84) and (623.22,136.14) .. (614.17,160.1) .. controls (605.12,184.05) and (591.71,201.59) .. (584.21,199.28) -- (363.95,131.22) .. controls (356.45,128.9) and (357.7,107.6) .. (366.75,83.65) .. controls (375.8,59.69) and (389.21,42.15) .. (396.71,44.47) .. controls (404.21,46.79) and (402.96,68.09) .. (393.91,92.04) .. controls (384.86,116) and (371.45,133.54) .. (363.95,131.22) ;
\draw  [color={rgb, 255:red, 139; green, 87; blue, 42 }  ,draw opacity=1 ][line width=1.5]  (457,160.94) .. controls (447.5,155.74) and (458.28,135.67) .. (481.07,116.12) .. controls (503.86,96.56) and (530.04,84.93) .. (539.54,90.13) .. controls (549.04,95.33) and (538.26,115.39) .. (515.47,134.95) .. controls (492.67,154.51) and (466.5,166.14) .. (457,160.94) -- cycle ;
\draw  [color={rgb, 255:red, 126; green, 211; blue, 33 }  ,draw opacity=1 ][line width=1.5]  (542.5,31.11) -- (542.5,266.66) .. controls (542.5,274.16) and (522.24,280.23) .. (497.25,280.23) .. controls (472.26,280.23) and (452,274.16) .. (452,266.66) -- (452,31.11) .. controls (452,23.62) and (472.26,17.54) .. (497.25,17.54) .. controls (522.24,17.54) and (542.5,23.62) .. (542.5,31.11) .. controls (542.5,38.61) and (522.24,44.69) .. (497.25,44.69) .. controls (472.26,44.69) and (452,38.61) .. (452,31.11) ;
\draw  [draw opacity=0][line width=1.5]  (542.1,264.38) .. controls (542.34,264.95) and (542.48,265.53) .. (542.49,266.13) .. controls (542.61,273.7) and (522.61,280.16) .. (497.83,280.56) .. controls (473.05,280.95) and (452.86,275.13) .. (452.74,267.55) .. controls (452.73,266.85) and (452.89,266.16) .. (453.21,265.49) -- (497.61,266.84) -- cycle ; \draw  [color={rgb, 255:red, 65; green, 117; blue, 5 }  ,draw opacity=1 ][line width=1.5]  (542.1,264.38) .. controls (542.34,264.95) and (542.48,265.53) .. (542.49,266.13) .. controls (542.61,273.7) and (522.61,280.16) .. (497.83,280.56) .. controls (473.05,280.95) and (452.86,275.13) .. (452.74,267.55) .. controls (452.73,266.85) and (452.89,266.16) .. (453.21,265.49) ;  
\draw  [draw opacity=0][dash pattern={on 5.63pt off 4.5pt}][line width=1.5]  (453.61,265.95) .. controls (453.37,265.38) and (453.25,264.79) .. (453.25,264.2) .. controls (453.29,256.62) and (473.41,250.57) .. (498.19,250.68) .. controls (522.98,250.79) and (543.04,257.03) .. (543.01,264.61) .. controls (543.01,265.31) and (542.83,265.99) .. (542.5,266.66) -- (498.13,264.4) -- cycle ; \draw  [color={rgb, 255:red, 65; green, 117; blue, 5 }  ,draw opacity=1 ][dash pattern={on 5.63pt off 4.5pt}][line width=1.5]  (453.61,265.95) .. controls (453.37,265.38) and (453.25,264.79) .. (453.25,264.2) .. controls (453.29,256.62) and (473.41,250.57) .. (498.19,250.68) .. controls (522.98,250.79) and (543.04,257.03) .. (543.01,264.61) .. controls (543.01,265.31) and (542.83,265.99) .. (542.5,266.66) ;  
\draw    (284.17,245.47) -- (218.03,355.69) ;
\draw [shift={(217,357.4)}, rotate = 300.97] [color={rgb, 255:red, 0; green, 0; blue, 0 }  ][line width=0.75]    (10.93,-3.29) .. controls (6.95,-1.4) and (3.31,-0.3) .. (0,0) .. controls (3.31,0.3) and (6.95,1.4) .. (10.93,3.29)   ;
\draw    (680,256.73) -- (735.02,353.99) ;
\draw [shift={(736,355.73)}, rotate = 240.51] [color={rgb, 255:red, 0; green, 0; blue, 0 }  ][line width=0.75]    (10.93,-3.29) .. controls (6.95,-1.4) and (3.31,-0.3) .. (0,0) .. controls (3.31,0.3) and (6.95,1.4) .. (10.93,3.29)   ;
\draw  [dash pattern={on 4.5pt off 4.5pt}]  (409.67,509.78) -- (588.83,509.92) ;
\draw [shift={(590.83,509.92)}, rotate = 180.04] [color={rgb, 255:red, 0; green, 0; blue, 0 }  ][line width=0.75]    (10.93,-3.29) .. controls (6.95,-1.4) and (3.31,-0.3) .. (0,0) .. controls (3.31,0.3) and (6.95,1.4) .. (10.93,3.29)   ;
\draw  [color={rgb, 255:red, 208; green, 2; blue, 27 }  ,draw opacity=1 ][line width=1.5]  (143.9,574.6) -- (346.77,574.6) -- (248.28,642.82) -- (45.41,642.82) -- cycle ;
\draw  [color={rgb, 255:red, 126; green, 211; blue, 33 }  ,draw opacity=1 ][line width=1.5]  (248.74,390.11) -- (248.74,608.06) .. controls (248.74,614.86) and (230.35,620.38) .. (207.68,620.38) .. controls (185,620.38) and (166.62,614.86) .. (166.62,608.06) -- (166.62,390.11) .. controls (166.62,383.31) and (185,377.8) .. (207.68,377.8) .. controls (230.35,377.8) and (248.74,383.31) .. (248.74,390.11) .. controls (248.74,396.92) and (230.35,402.43) .. (207.68,402.43) .. controls (185,402.43) and (166.62,396.92) .. (166.62,390.11) ;
\draw  [draw opacity=0][line width=1.5]  (248.38,605.76) .. controls (248.6,606.28) and (248.71,606.81) .. (248.72,607.34) .. controls (248.83,614.34) and (230.69,620.31) .. (208.21,620.67) .. controls (185.72,621.04) and (167.4,615.66) .. (167.3,608.67) .. controls (167.29,608.03) and (167.43,607.41) .. (167.7,606.79) -- (208.01,608.01) -- cycle ; \draw  [color={rgb, 255:red, 65; green, 117; blue, 5 }  ,draw opacity=1 ][line width=1.5]  (248.38,605.76) .. controls (248.6,606.28) and (248.71,606.81) .. (248.72,607.34) .. controls (248.83,614.34) and (230.69,620.31) .. (208.21,620.67) .. controls (185.72,621.04) and (167.4,615.66) .. (167.3,608.67) .. controls (167.29,608.03) and (167.43,607.41) .. (167.7,606.79) ;  
\draw  [draw opacity=0][dash pattern={on 5.63pt off 4.5pt}][line width=1.5]  (168.07,607.16) .. controls (167.86,606.64) and (167.76,606.11) .. (167.76,605.57) .. controls (167.79,598.58) and (186.05,592.99) .. (208.54,593.09) .. controls (231.02,593.19) and (249.23,598.95) .. (249.2,605.94) .. controls (249.19,606.58) and (249.04,607.2) .. (248.75,607.8) -- (208.48,605.76) -- cycle ; \draw  [color={rgb, 255:red, 65; green, 117; blue, 5 }  ,draw opacity=1 ][dash pattern={on 5.63pt off 4.5pt}][line width=1.5]  (168.07,607.16) .. controls (167.86,606.64) and (167.76,606.11) .. (167.76,605.57) .. controls (167.79,598.58) and (186.05,592.99) .. (208.54,593.09) .. controls (231.02,593.19) and (249.23,598.95) .. (249.2,605.94) .. controls (249.19,606.58) and (249.04,607.2) .. (248.75,607.8) ;  
\draw  [draw opacity=0][line width=1.5]  (249.25,476.4) .. controls (249.31,476.66) and (249.34,476.92) .. (249.34,477.18) .. controls (249.33,483.44) and (230.94,488.32) .. (208.27,488.09) .. controls (185.87,487.86) and (167.66,482.71) .. (167.26,476.55) -- (208.3,476.75) -- cycle ; \draw  [color={rgb, 255:red, 144; green, 19; blue, 254 }  ,draw opacity=1 ][line width=1.5]  (249.25,476.4) .. controls (249.31,476.66) and (249.34,476.92) .. (249.34,477.18) .. controls (249.33,483.44) and (230.94,488.32) .. (208.27,488.09) .. controls (185.87,487.86) and (167.66,482.71) .. (167.26,476.55) ;  
\draw  [draw opacity=0][dash pattern={on 5.63pt off 4.5pt}][line width=1.5]  (167.47,474.54) .. controls (169.2,467.01) and (186.94,461.51) .. (208.5,462.08) .. controls (228.88,462.61) and (245.72,468.39) .. (248.77,475.46) -- (208.3,476.75) -- cycle ; \draw  [color={rgb, 255:red, 144; green, 19; blue, 254 }  ,draw opacity=1 ][dash pattern={on 5.63pt off 4.5pt}][line width=1.5]  (167.47,474.54) .. controls (169.2,467.01) and (186.94,461.51) .. (208.5,462.08) .. controls (228.88,462.61) and (245.72,468.39) .. (248.77,475.46) ;  
\draw  [color={rgb, 255:red, 208; green, 2; blue, 27 }  ,draw opacity=1 ][line width=1.5]  (719.23,479.92) -- (942.82,479.92) -- (834.26,543.61) -- (610.67,543.61) -- cycle ;
\draw  [color={rgb, 255:red, 245; green, 166; blue, 35 }  ,draw opacity=1 ][line width=1.5]  (801.16,418.22) -- (801.16,597.85) .. controls (801.16,604.53) and (783.47,609.94) .. (761.66,609.94) .. controls (739.84,609.94) and (722.16,604.53) .. (722.16,597.85) -- (722.16,418.22) .. controls (722.16,411.55) and (739.84,406.13) .. (761.66,406.13) .. controls (783.47,406.13) and (801.16,411.55) .. (801.16,418.22) .. controls (801.16,424.9) and (783.47,430.31) .. (761.66,430.31) .. controls (739.84,430.31) and (722.16,424.9) .. (722.16,418.22) ;
\draw  [color={rgb, 255:red, 126; green, 211; blue, 33 }  ,draw opacity=1 ][line width=1.5]  (722.55,511.3) .. controls (722.55,503.5) and (740.24,497.18) .. (762.05,497.18) .. controls (783.87,497.18) and (801.55,503.5) .. (801.55,511.3) .. controls (801.55,519.1) and (783.87,525.43) .. (762.05,525.43) .. controls (740.24,525.43) and (722.55,519.1) .. (722.55,511.3) -- cycle ;
\draw [color={rgb, 255:red, 245; green, 166; blue, 35 }  ,draw opacity=1 ][line width=1.5]    (166.31,492.42) .. controls (166.6,481.91) and (182.07,455.2) .. (194.87,453.6) .. controls (207.67,452) and (189.69,499.38) .. (199.67,496) .. controls (209.65,492.61) and (208.75,469.54) .. (216.47,464) .. controls (224.18,458.46) and (222.22,504.08) .. (227.72,503.77) .. controls (233.22,503.47) and (245.27,428.8) .. (248.77,475.46) .. controls (252.27,522.12) and (234.87,522.4) .. (207.67,518.4) .. controls (180.47,514.4) and (167.25,499.71) .. (166.31,492.42) -- cycle ;

\draw (633.48,199.06) node [anchor=north west][inner sep=0.75pt]  [font=\Large,color={rgb, 255:red, 208; green, 2; blue, 27 }  ,opacity=1 ]  {$E_{P} \simeq \mathbb{P}^{2}$};
\draw (558.15,17.72) node [anchor=north west][inner sep=0.75pt]  [font=\Large,color={rgb, 255:red, 126; green, 211; blue, 33 }  ,opacity=1 ]  {$E_{e} \simeq \mathbb{F}_{6}$};
\draw (555.15,68.91) node [anchor=north west][inner sep=0.75pt]  [font=\Large,color={rgb, 255:red, 139; green, 87; blue, 42 }  ,opacity=1 ]  {$E_{e} \ \cap \ E_{f}$};
\draw (223.2,278.87) node [anchor=north west][inner sep=0.75pt]  [font=\Large]  {$\sigma _{f}$};
\draw (725.33,298.4) node [anchor=north west][inner sep=0.75pt]  [font=\Large]  {$\beta $};
\draw (334.47,587.31) node [anchor=north west][inner sep=0.75pt]  [font=\Large,color={rgb, 255:red, 208; green, 2; blue, 27 }  ,opacity=1 ]  {$E_{P} \simeq \mathbb{P}^{2}$};
\draw (57.24,481.39) node [anchor=north west][inner sep=0.75pt]  [font=\Large,color={rgb, 255:red, 245; green, 166; blue, 35 }  ,opacity=1 ]  {$f\ =\ \alpha ( e')$};
\draw (262.02,376.7) node [anchor=north west][inner sep=0.75pt]  [font=\Large,color={rgb, 255:red, 126; green, 211; blue, 33 }  ,opacity=1 ]  {$E_{e} \simeq \mathbb{F}_{6}$};
\draw (263.6,527.2) node [anchor=north west][inner sep=0.75pt]  [font=\Large,color={rgb, 255:red, 126; green, 211; blue, 33 }  ,opacity=1 ]  {$\{A=0\}$};
\draw (278.82,53.06) node [anchor=north west][inner sep=0.75pt]  [font=\Large,color={rgb, 255:red, 245; green, 166; blue, 35 }  ,opacity=1 ]  {$E_{f} \simeq \mathbb{F}_{2}$};
\draw (642.04,383.77) node [anchor=north west][inner sep=0.75pt]  [font=\Large,color={rgb, 255:red, 245; green, 166; blue, 35 }  ,opacity=1 ]  {$E_{f} \simeq \mathbb{F}_{2}$};
\draw (896.28,443.01) node [anchor=north west][inner sep=0.75pt]  [font=\Large,color={rgb, 255:red, 208; green, 2; blue, 27 }  ,opacity=1 ]  {$E_{P} \simeq \mathbb{P}^{2}$};
\draw (802.81,482.71) node [anchor=north west][inner sep=0.75pt]  [font=\Large,color={rgb, 255:red, 126; green, 211; blue, 33 }  ,opacity=1 ]  {$E_{f} \ \cap \ E_{P}$};
\draw (370.6,364.73) node [anchor=north west][inner sep=0.75pt]  [font=\Large]  {$Z_{2}$};
\draw (912.6,365.73) node [anchor=north west][inner sep=0.75pt]  [font=\Large]  {$Z_{3}$};
\draw (679.6,18.81) node [anchor=north west][inner sep=0.75pt]  [font=\Large]  {$Y_{2}$};

\end{tikzpicture}

}  
\caption{Third link in a Sarkisov factorization of $\phi$.}
\label{fig:third link}
\end{figure}

\subparagraph{Stage 4:} The computation of the Sarkisov degree of the induced birational map $\phi \circ \sigma_P \circ \sigma_e \circ \alpha^{-1} \circ \sigma_f \circ \beta^{-1}$ will be a little different. The issue here is that the pair $(Z_3,c\h_{Z_3})$ is canonical for any positive rational number $c$, since $\h_{Z_3}$ is base point free. So the notion of canonical threshold would lead us to $c_3~ \text{``$=$''}~\infty$. This is in accordance with Definition \ref{defn sark deg dim 3}.

In this case, the Sarkisov degree $(\mu_3,c_3,e_3)$ becomes $\left( \dfrac{1}{2},\infty,* \right)$ and it is smaller than $(\mu_2,c_2,e_2)$. 

Since $c_3= \infty\  \text{``$\geq$''}\  2=\dfrac{1}{\mu_3}$, the fourth link in the Sarkisov factorization is of type III or IV. We need to run the $(K_{Z_3}+2\h_{Z_3})$-MMP over $\Spec(\C)$. This log MMP results exactly in the divisorial contraction given by the blowup of $\p^3$ at $P$, that is, the map $\sigma_P$. We observe that $E_f$ is precisely $\Exc(\sigma_P)$.  

One can compute that the Sarkisov degree $(\mu_4,c_4,e_4)$ of the induced birational map $\phi \circ \sigma_P \circ \sigma_e \circ \alpha^{-1} \circ \sigma_f \circ \beta^{-1} \circ \sigma_P^{-1}$ is $\left( \dfrac{1}{4},\infty,* \right)$, which is smaller than the previous one. Such Sarkisov degree implies that this map is an automorphism of $\p^3$.

The (volume preserving) factorization of $\phi$ is ended up by the Sarkisov link of type III given by $\sigma_P$. We have the following picture:

\vspace{20.0cm}

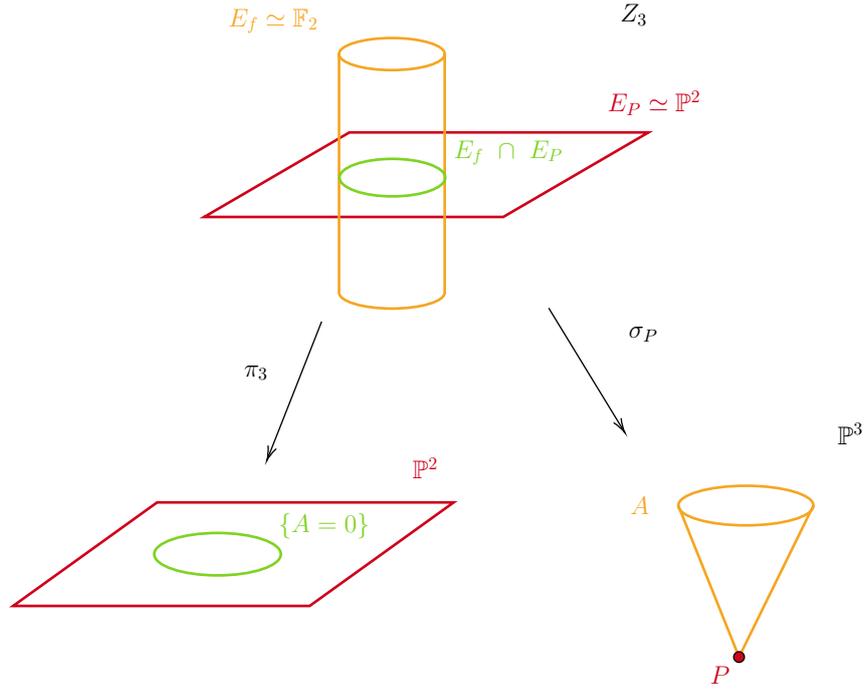
\begin{figure}[htb]
\centering
\resizebox{0.7\textwidth}{!}%
{

\tikzset{every picture/.style={line width=0.75pt}} 

\begin{tikzpicture}[x=0.75pt,y=0.75pt,yscale=-1,xscale=1]

\draw  [color={rgb, 255:red, 245; green, 166; blue, 35 }  ,draw opacity=1 ][line width=1.5]  (528.13,404.98) .. controls (528.13,396.74) and (550.74,390.07) .. (578.63,390.07) .. controls (606.52,390.07) and (629.13,396.74) .. (629.13,404.98) .. controls (629.13,413.21) and (606.52,419.89) .. (578.63,419.89) .. controls (550.74,419.89) and (528.13,413.21) .. (528.13,404.98) -- cycle ;
\draw [color={rgb, 255:red, 245; green, 166; blue, 35 }  ,draw opacity=1 ][line width=1.5]    (528.13,404.98) -- (573.52,519.07) -- (629.13,404.98) ;
\draw  [color={rgb, 255:red, 208; green, 2; blue, 27 }  ,draw opacity=1 ][line width=1.5]  (138.48,402.6) -- (360.18,402.6) -- (252.53,480.6) -- (30.83,480.6) -- cycle ;
\draw  [color={rgb, 255:red, 126; green, 211; blue, 33 }  ,draw opacity=1 ][line width=1.5]  (136.33,441.6) .. controls (136.33,432.77) and (157.53,425.62) .. (183.67,425.62) .. controls (209.81,425.62) and (231,432.77) .. (231,441.6) .. controls (231,450.43) and (209.81,457.58) .. (183.67,457.58) .. controls (157.53,457.58) and (136.33,450.43) .. (136.33,441.6) -- cycle ;
\draw    (431.23,256.17) -- (487.2,348.46) ;
\draw [shift={(488.23,350.17)}, rotate = 238.77] [color={rgb, 255:red, 0; green, 0; blue, 0 }  ][line width=0.75]    (10.93,-3.29) .. controls (6.95,-1.4) and (3.31,-0.3) .. (0,0) .. controls (3.31,0.3) and (6.95,1.4) .. (10.93,3.29)   ;
\draw    (261.5,266.35) -- (221.23,369.26) ;
\draw [shift={(220.5,371.13)}, rotate = 291.37] [color={rgb, 255:red, 0; green, 0; blue, 0 }  ][line width=0.75]    (10.93,-3.29) .. controls (6.95,-1.4) and (3.31,-0.3) .. (0,0) .. controls (3.31,0.3) and (6.95,1.4) .. (10.93,3.29)   ;
\draw  [fill={rgb, 255:red, 208; green, 2; blue, 27 }  ,fill opacity=1 ] (569.52,519.07) .. controls (569.52,516.86) and (571.31,515.07) .. (573.52,515.07) .. controls (575.73,515.07) and (577.52,516.86) .. (577.52,519.07) .. controls (577.52,521.28) and (575.73,523.07) .. (573.52,523.07) .. controls (571.31,523.07) and (569.52,521.28) .. (569.52,519.07) -- cycle ;
\draw  [color={rgb, 255:red, 208; green, 2; blue, 27 }  ,draw opacity=1 ][line width=1.5]  (282.23,123.88) -- (505.82,123.88) -- (397.26,187.57) -- (173.67,187.57) -- cycle ;
\draw  [color={rgb, 255:red, 245; green, 166; blue, 35 }  ,draw opacity=1 ][line width=1.5]  (353.49,64.85) -- (353.49,244.48) .. controls (353.49,251.16) and (335.81,256.57) .. (313.99,256.57) .. controls (292.17,256.57) and (274.49,251.16) .. (274.49,244.48) -- (274.49,64.85) .. controls (274.49,58.17) and (292.17,52.76) .. (313.99,52.76) .. controls (335.81,52.76) and (353.49,58.17) .. (353.49,64.85) .. controls (353.49,71.53) and (335.81,76.94) .. (313.99,76.94) .. controls (292.17,76.94) and (274.49,71.53) .. (274.49,64.85) ;
\draw  [color={rgb, 255:red, 126; green, 211; blue, 33 }  ,draw opacity=1 ][line width=1.5]  (274.88,157.93) .. controls (274.88,150.13) and (292.57,143.81) .. (314.38,143.81) .. controls (336.2,143.81) and (353.89,150.13) .. (353.89,157.93) .. controls (353.89,165.73) and (336.2,172.06) .. (314.38,172.06) .. controls (292.57,172.06) and (274.88,165.73) .. (274.88,157.93) -- cycle ;

\draw (550.53,524.75) node [anchor=north west][inner sep=0.75pt]  [font=\Large,color={rgb, 255:red, 208; green, 2; blue, 27 }  ,opacity=1 ]  {$P$};
\draw (490.8,396.75) node [anchor=north west][inner sep=0.75pt]  [font=\Large,color={rgb, 255:red, 245; green, 166; blue, 35 }  ,opacity=1 ]  {$A$};
\draw (227.93,409.2) node [anchor=north west][inner sep=0.75pt]  [font=\Large,color={rgb, 255:red, 126; green, 211; blue, 33 }  ,opacity=1 ]  {$\{A=0\}$};
\draw (328.2,367.67) node [anchor=north west][inner sep=0.75pt]  [font=\Large,color={rgb, 255:red, 208; green, 2; blue, 27 }  ,opacity=1 ]  {$\mathbb{P}^{2}$};
\draw (489.6,268.53) node [anchor=north west][inner sep=0.75pt]  [font=\Large]  {$\sigma _{P}$};
\draw (202,297.73) node [anchor=north west][inner sep=0.75pt]  [font=\Large]  {$\pi _{3}$};
\draw (190.37,29.06) node [anchor=north west][inner sep=0.75pt]  [font=\Large,color={rgb, 255:red, 245; green, 166; blue, 35 }  ,opacity=1 ]  {$E_{f} \simeq \mathbb{F}_{2}$};
\draw (473.94,90.97) node [anchor=north west][inner sep=0.75pt]  [font=\Large,color={rgb, 255:red, 208; green, 2; blue, 27 }  ,opacity=1 ]  {$E_{P} \simeq \mathbb{P}^{2}$};
\draw (358.81,128.33) node [anchor=north west][inner sep=0.75pt]  [font=\Large,color={rgb, 255:red, 126; green, 211; blue, 33 }  ,opacity=1 ]  {$E_{f} \ \cap \ E_{P}$};
\draw (483.27,25.73) node [anchor=north west][inner sep=0.75pt]  [font=\Large]  {$Z_{3}$};
\draw (645.93,341.73) node [anchor=north west][inner sep=0.75pt]  [font=\Large]  {$\mathbb{P}^{3}$};

\end{tikzpicture}
}  
\caption{Fourth link in a Sarkisov factorization of $\phi$.}
\label{fig:fourth link}
\end{figure}

We observe that the strict transforms of $D$ remained nonsingular and isomorphic to $D$ along the intermediate stages. In terms of a commutative diagram, this volume preserving factorization of $\phi$ is expressed in the following way:
\[\begin{tikzcd}[column sep=0.005cm]
	&& {Y_1} && {Y_2} \\
	& {Z_1 \simeq \mathbb{P}(\mathcal{O}_{\mathbb{P}^2} \oplus \mathcal{O}_{\mathbb{P}^2}(1))} && {Z_2 \simeq \mathbb{P}(\mathcal{O}_{\mathbb{P}^2} \oplus \mathcal{O}_{\mathbb{P}^2}(3))} && {Z_3 \simeq \mathbb{P}(\mathcal{O}_{\mathbb{P}^2} \oplus \mathcal{O}_{\mathbb{P}^2}(1))} \\
	{\mathbb{P}^3} & {\mathbb{P}^2} && {\mathbb{P}^2} && {\mathbb{P}^2} & {\mathbb{P}^3} \\
	{\text{Spec}(\mathbb{C})} &&&&&& {\text{Spec}(\mathbb{C})}
	\arrow[dashed, from=2-2, to=2-4]
	\arrow["{\sigma_e}"', from=1-3, to=2-2]
	\arrow["\alpha", from=1-3, to=2-4]
	\arrow["{\sigma_f}"', from=1-5, to=2-4]
	\arrow["\beta", from=1-5, to=2-6]
	\arrow[dashed, from=2-4, to=2-6]
	\arrow["{\pi_1}"', from=2-2, to=3-2]
	\arrow["{\pi_2}"', from=2-4, to=3-4]
	\arrow["{\pi_3}"', from=2-6, to=3-6]
	\arrow[Rightarrow, no head, from=3-2, to=3-4]
	\arrow[Rightarrow, no head, from=3-4, to=3-6]
	\arrow["{\sigma_P}"', from=2-2, to=3-1]
	\arrow["{\sigma_P}", from=2-6, to=3-7]
	\arrow["\phi", curve={height=30pt}, dashed, from=3-1, to=3-7]
	\arrow[from=3-7, to=4-7]
	\arrow[from=3-1, to=4-1]
\end{tikzcd}\]

We have the following sequence of Sarkisov degrees of the induced birational maps:
\begin{center}
$\left( \dfrac{3}{4},1,9 \right) > \left( \dfrac{1}{2},1,8 \right) > \left( \dfrac{1}{2},1,1 \right) > \left( \dfrac{1}{2},\infty,* \right) > \left( \dfrac{1}{4},\infty,* \right)$.
\end{center}

\paragraph{Non-volume preserving factorization of $\phi$.} Let us make some comments if instead of proceeding with $\sigma_e$ in Stage 2, we had chosen $\sigma_1$. We invite the reader to fill out the details and make drawings in order to see what is happening geometrically.

For $i \in \{1,\ldots,6\} $, denote $P_i \coloneqq L_i \cap E_P$. In this other Sarkisov factorization the $(i+1)$-th Sarkisov link is of type I, and it is given by $\sigma_i \colon \Bl_{L_1,\ldots,L_i}(Z_1) \rightarrow \Bl_{L_1,\ldots,L_{i-1}}(Z_1)$, where $\Bl_{L_1,\ldots,L_i}(Z_1) \rightarrow \Bl_{P_1,\ldots,P_i}(\p^2)$ is the corresponding structure of Mori fibered space. One can show that $E_i \coloneqq \Exc(\sigma_i)$ is isomorphic to $\F_0 \simeq \p^1 \times \p^1$ for all $i \in \{1,\ldots,6\}$. 

We have the following picture for $\sigma_1$, where $g=\tilde{E_P} \cap E_1 \simeq \p^1$. Geometrically, $g$ represents all the normal directions to $L_1$ at $P_1$.

\vspace{20.0cm}

\begin{figure}[htb]
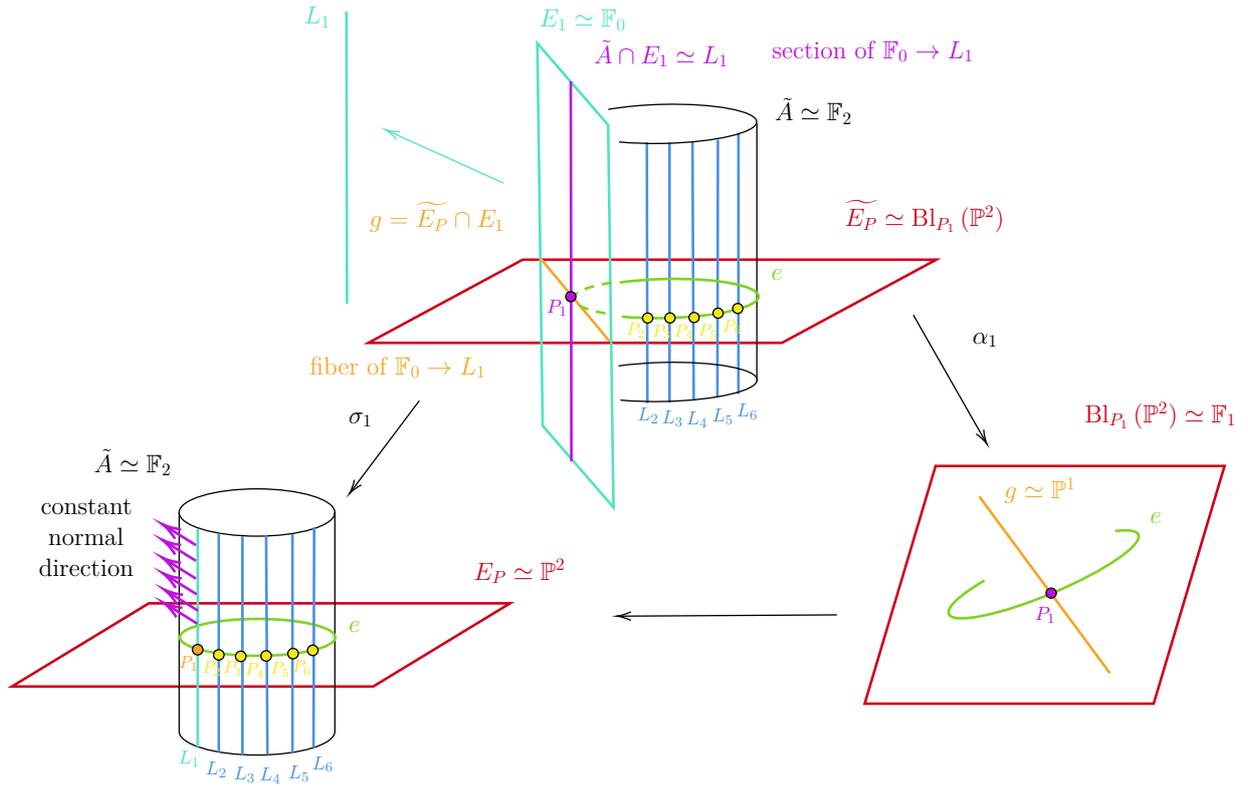

\centering
\resizebox{1.0\textwidth}{!}%
{

\tikzset{every picture/.style={line width=0.75pt}} 



}  
\caption{Sarkisov link $\sigma_1$.}
\label{fig:second link non vp}
\end{figure}

After $\sigma_6$, the next two Sarkisov links are analogous to the intermediate ones in the volume preserving factorization described previously. In particular, they are volume preserving Sarkisov links of type II. 

\begin{figure}[htb]
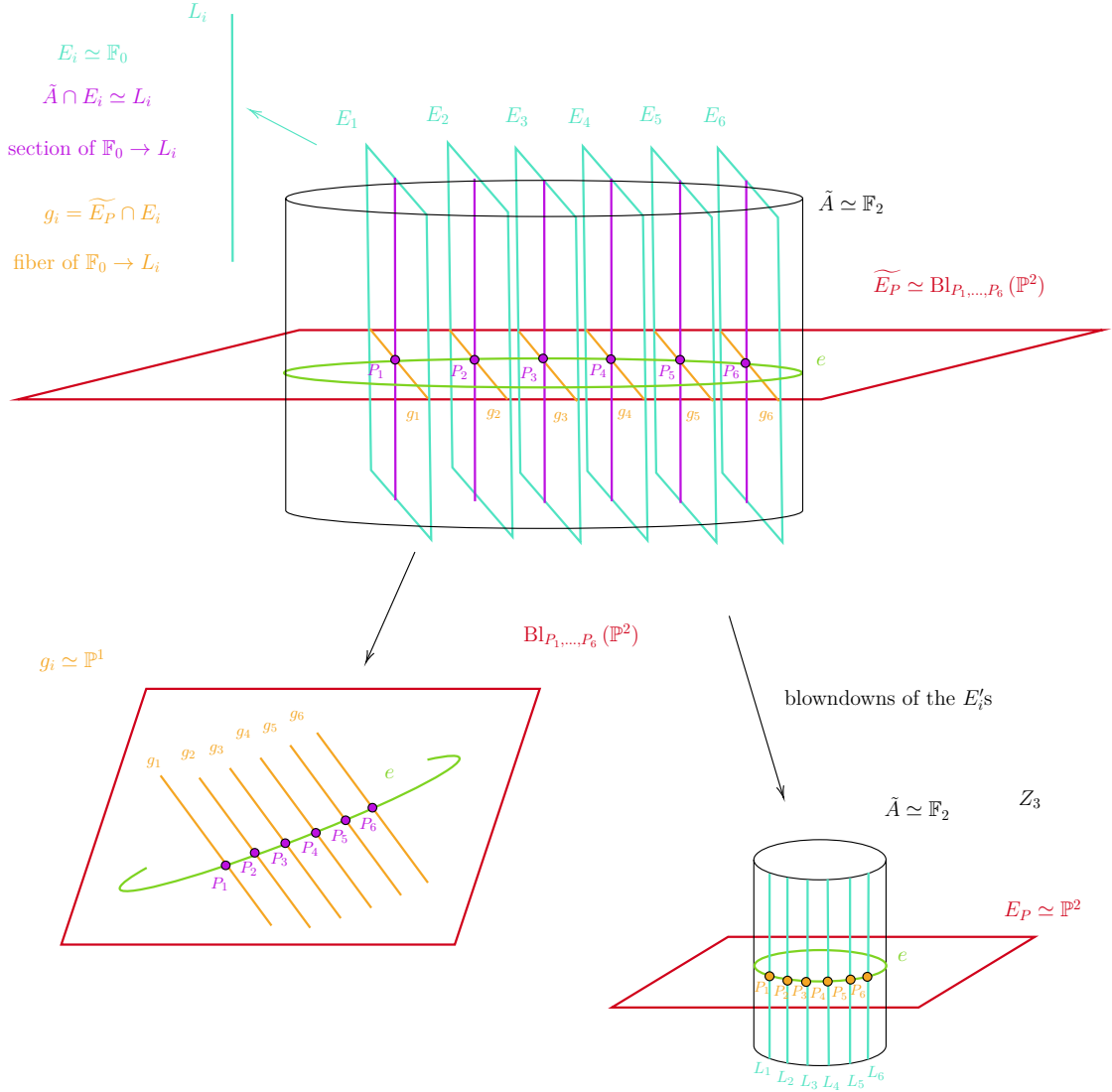

\centering
\resizebox{0.9\textwidth}{!}%
{

\tikzset{every picture/.style={line width=0.75pt}} 



}  
\caption{Blowdowns of the $E_i$'s.}
\label{fig:blowdown of the E_i's}
\end{figure}

Finally, the remaining ones are given by the consecutive blowdowns of the corresponding strict transforms of $E_1, \ldots, E_6$ and $E_P$, respectively. All of them are Sarkisov links of type III and the last one is volume preserving. The setting is depicted in Figure \ref{fig:blowdown of the E_i's}.

\vspace{20.0cm}

We have the following sequence of Sarkisov degrees of the induced birational maps:

\begin{center}
    $\left( \dfrac{3}{4},1,9 \right) > \underbrace{\left( \dfrac{1}{2},1,8 \right) > \ldots > \left( \dfrac{1}{2},1,2 \right)}_{\text{relative to the $\sigma_i$'s}} > \left( \dfrac{1}{2},1,1 \right) > \underbrace{\left( \dfrac{1}{2},2,6 \right) > \ldots > \left( \dfrac{1}{2},2,1 \right)}_{\text{relative to the blowdowns of the $E_i$'s}}$\\
    $>\left( \dfrac{1}{2},\infty,* \right) > \left( \dfrac{1}{4},\infty,* \right)$.
\end{center}

We remark that it is relevant that $P \in \Sing(\h)$ is a singularity of type $A_1$ for the increasing of the canonical threshold in Stage 8. The induced birational map before the Sarkisov links of type III has a base locus given by $\{P_1,\ldots,P_6 \}$.

\subsection{Conclusion} By the previous detailed example, we can see the existence of a Sarkisov factorization for a volume preserving map that is not automatically volume preserving. This means that Theorem  \ref{thm sark=vp cubic} is very particular for dimension 2. 

The point is that the volume preserving factorization assured by \cite[Theorem 1.1]{ck} is induced by a standard Sarkisov factorization constructed in a special way \cite[Theorem 3.3]{ck}. But this special factorization may not be obtained by following algorithmically the Sarkisov program, at least in dimension 3, if we make certain choices along the process.

In a careful analysis of the volume preserving decomposition of $\phi$ exhibited in \cite{acm}, we show that it can be obtained by choosing conveniently the centers of the extremal blowups initiating Sarkisov links of type I and II. Furthermore, this factorization has the effect of resolving the map $\phi$ along the process as a consequence of the untwisting of the Sarkisov degree of the induced birational maps. 

\let\oldbibliography\thebibliography
\renewcommand{\thebibliography}[1]{\oldbibliography{#1}
\setlength{\itemsep}{0pt}} 

\end{document}